\documentclass[natbib]{svjour3}
\smartqed  
\usepackage{graphicx,natbib}
\usepackage{amsfonts,amsmath,amssymb,mathrsfs,epsfig,a4}
\usepackage{latexsym}
\usepackage{xcolor}
\usepackage{longtable}

\usepackage{amsthm}
\usepackage{mathabx}
\usepackage{multirow}

\allowdisplaybreaks

\newcommand{\bN}{\mathbb{N}}
\newcommand{\bP}{\mathbb{P}}

\newcommand{\bR}{\mathbb{R}}

\newcommand{\bH}{\mathbb{H}}

\newcommand{\cA}{\mathcal{A}}
\newcommand{\cB}{\mathcal{B}}

\newcommand{\cM}{\mathcal{M}}

\newcommand{\cS}{\mathcal{S}}

\newcommand{\Sd}{\mathbb{S}_{d-1}}
\newcommand{\Seins}{\mathbb{S}_{1}}

\newcommand{\CoMed}{\text{\sc MMed}}

\newcommand{\QMed}{\text{\sc QMed}}
\newcommand{\TMed}{\text{\sc TMed}}
\newcommand{\OMed}{\text{\sc OMed}}
\newcommand{\SMed}{\text{\sc SMed}}
\newcommand{\SO}{\mathbb{SO}}
\newcommand{\OO}{\mathbb{O}}
\newcommand{\so}{\mathfrak{so}}
\newcommand{\ML}{\text{\sc ML}}

\newcommand{\unif}{\text{\rm unif}}
\newcommand{\Med}{\text{\sc Med}}
\newcommand{\mMed}{\text{\,\sc MMed}}
\newcommand{\cov}{\text{\rm cov}}
\newcommand{\diag}{\text{\rm diag}}
\newcommand{\Ell}{\text{\rm Ell}}
\newcommand{\Sym}{\text{\rm Sym}}

\renewcommand{\t}{^{\text{t}}}

\newcommand{\dconv}{\to_{\text{\tiny d}}}

\journalname{Metrika}

\makeatletter
\def\LT@makecaption#1#2#3{%
  \LT@mcol\LT@cols c{\hbox to\z@{\hss\parbox[t]\LTcapwidth{%
    \captionstyle
    \sbox\@tempboxa{{\floatlegendstyle#1{#2}\floatcounterend}\capstrut #3}%
    \ifdim\wd\@tempboxa>\hsize
      {\floatlegendstyle#1{#2}\floatcounterend}\capstrut #3\par%
    \else
      \hbox to\hsize{\leftlegendglue\box\@tempboxa\hfil}%
    \fi
    \endgraf\vskip\baselineskip}%
  \hss}}}
\LTcapwidth=\dimexpr\textwidth-2\fboxsep\relax
\makeatother

\begin{document}

\title{The quarter median}
	
\author{Ludwig Baringhaus and Rudolf Gr\"ubel}

\institute{L. Baringhaus \at 
           Institute of Actuarial and Financial Mathematics  \\ 
	   Leibniz Universit\"at Hannover\\
	   Postfach 6009\\
	   30060 Hannover\\
	   Germany\\
           \email{lbaring@stochastik.uni-hannover.de} 
           \and
           R. Gr\"ubel \at
           Institute of Actuarial and Financial Mathematics \\ 
 	   Leibniz Universit\"at Hannover\\
 	   Postfach 6009\\
           D-30060 Hannover, Germany
 	   \email{rgrubel@stochastik.uni-hannover.de}
           }

\date{Received: date / Accepted: date}

\maketitle

\begin{abstract} 
We introduce and discuss a multivariate version of the classical median 
that is based on an equipartition property with respect to quarter spaces. 
These arise as pairwise intersections of the half-spaces associated with 
the coordinate hyperplanes of an orthogonal basis. We obtain results on 
existence, equivariance, and asymptotic normality.
\end{abstract}

\section{Introduction}\label{sec:intro} 

The classical one-dimensional median associated with a probability distribution 
$P$ on the Borel subsets of the real line $\bR$ is any value $\theta$ that satisfies the inequalities
\begin{equation}\label{eq:meddefone}
	P\bigl((-\infty,\theta]\bigr) \ge 1/2\,, \quad 
	P\bigl([\theta,\infty)\bigr) \ge 1/2\,.
\end{equation}
Under certain conditions, for example if $P$ has a positive density, we have equalities 
in \eqref{eq:meddefone} and the median is unique. In general, the set of solutions 
is the bounded closed interval $\left [\Med_-(P),\Med_+(P)\right ]$ with
\begin{align*}
  \Med_-(P)\,&:=\,\inf\bigl\{x\in\bR\!: P\bigl((-\infty,x]\bigr)\ge1/2\bigr\},\\
  \Med_+(P)\,&:=\,\sup\bigl\{x\in\bR\!: P\bigl([x,\infty)\bigr)\ge 1/2\bigr\}
\end{align*}
as the smallest and the largest median of $P$. Often, `the' 
median of $P$ is defined as the corresponding midpoint,
\begin{equation}\label{eq:midpoint1}
  \Med_m(P):= \bigl(\Med_-(P) + \Med_+(P)\bigr)/2.
\end{equation} 
For a data set $x_1,\ldots,x_n$ of real
numbers a (sample) median is any (distributional) median associated with the empirical 
distribution $P_n:=n^{-1} \sum_{k=1}^n \delta_{x_k}$, where $\delta_x$ denotes unit mass 
at $x$. With $x_{(1)}\le \ldots \le x_{(n)}$ as the values obtained by arranging
$x_1,\ldots,x_n$ in increasing order,
\begin{equation*}
  \Med_m(P_n)=\Med_-(P_n)=\Med_+(P_n)=x_{\left ((n+1)/2\right )}
\end{equation*}
is the unique sample median if $n$ is odd, and we get $\bigl[x_{(n/2)},x_{(n/2+1)}\bigr]$ 
as the set of sample medians if $n$ is even. 

Historically, the first attempt to generalize this to dimensions $d>1$ was to use coordinates:
If $x_1,\ldots,x_n$ is a subset of $\bR^d$ for some $d\ge 2$ then a coordinatewise  or 
\emph{marginal  median} of the empirical distribution $P_n:=n^{-1} \sum_{k=1}^n \delta_{x_k}$
is any vector whose $i$th coordinate is a (one-dimensional) median of the respective 
$i$th coordinates of the data vectors. This is the sample version; for
a probability measure $P$ on the Borel subsets $\cB^d$ of $\bR^d$ the push-forwards
$P^{\pi_i}$ of $P$ under the coordinate projections $\pi_i$, $i=1,\ldots,d$, would be used.
Formally, the marginal medians of $P$ are elements of the $d$-dimensional interval
$\bigl[\mMed_-(P),\mMed_+(P)\bigr]$ where $\mMed_\pm(P)$ have $\Med_\pm(P^{\pi_i})$ 
as their $i$th component, $i=1,\ldots,d$, and 
\begin{equation}\label{eq:midpointd}
     \mMed_m(P):=\bigl(\mMed_-(P)+\mMed_+(P)\bigr)/2
\end{equation}
is the midpoint marginal median.
It is well known that such component-wise extensions of the one-dimensional median are
not equivariant with respect to orthogonal transformations and thus depend on the coordinate
system chosen for the representation of the data.
One way to repair this is by `averaging out' (symmetrization), see \cite{Gr1996}. 

\begin{figure}
\vspace{-.5cm}
\setlength{\abovecaptionskip}{-0.5cm}
\includegraphics[scale=.6]{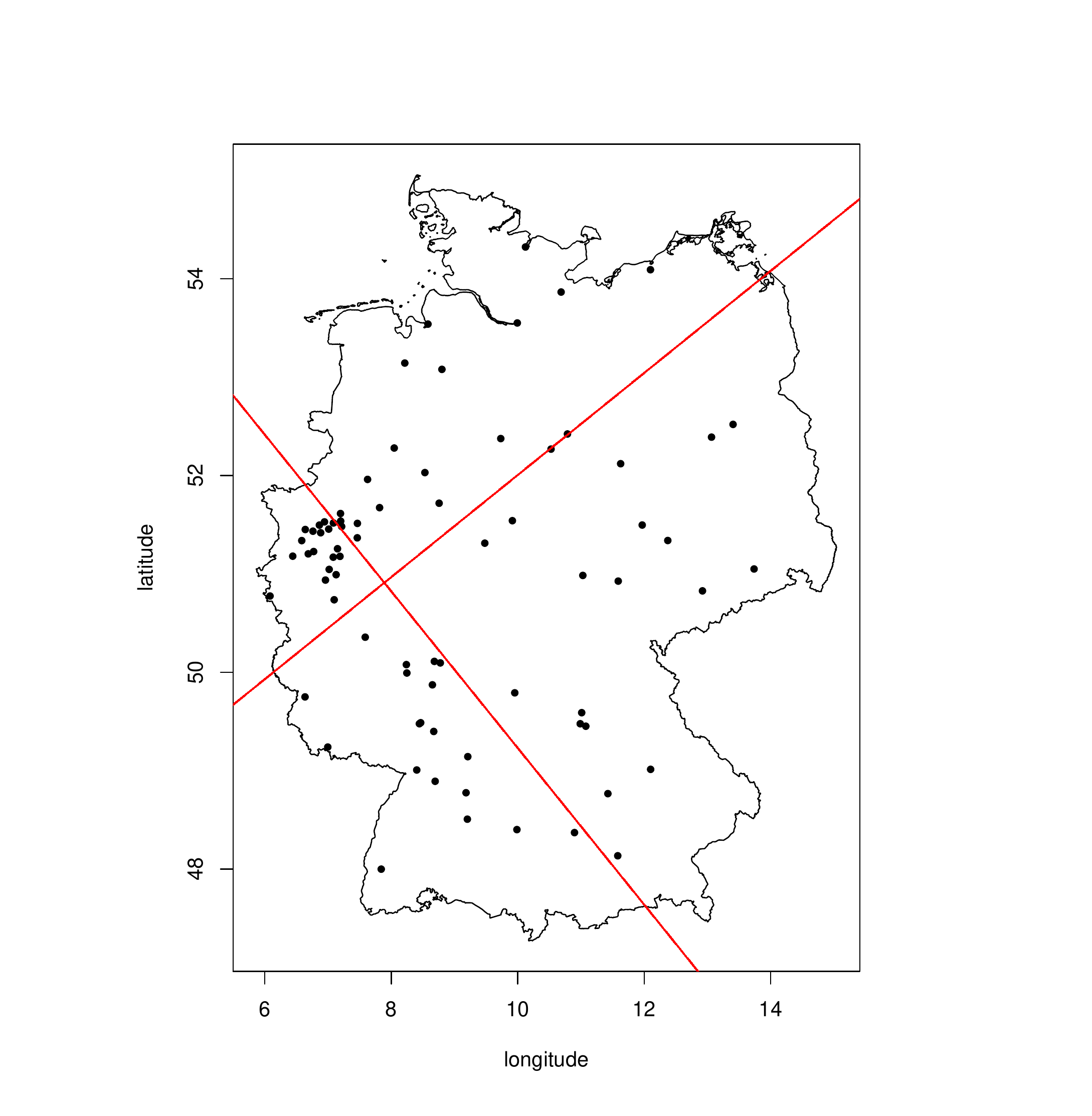}
\caption{The city data (see text for details)}
\label{fig:staedte}
\end{figure}

Other generalizations of the one-dimensional median to higher dimensions start  
with some characterizing property of the classical version 
and then use a multidimensional analogue. For example, it is known that in the 
one-dimensional case any median $\theta$ minimizes the function $x\mapsto \int |y-x|\, P(dy)$, 
which leads to the \emph{spatial} or \emph{$L^1$-median} of the data set, defined as a minimizer of
$x\mapsto \sum_{k=1}^n\| x-x_k\|$, where $\|\, \cdot\,\|$ denotes Euclidean distance. 
As this distance is invariant under orthogonal transformation, the spatial median does
not depend on the chosen coordinates. 

Closer to the above equipartition property are 
the various concepts related to data depth, where for each $x\in\bR^d$ and each affine hyperplane 
$H$ with $x\in H$ the number of data points to one side of $H$ is the basic ingredient: 
A centerpoint median is defined as any such $x$ with the property that the minimum number is 
at least a fraction $n/(d+1)$ of the data, and Tukey's  median is a point that maximizes these 
data depth values; see e.g.~\cite{DG1992}. This whole area has attracted the attention of 
many researchers, often with emphasis on robustness. The literature on 
multivariate medians up to about 1990 is covered in the review of~\cite{Sm1990}, which 
includes an interesting discussion of the history of the concept.
For more recent reviews of this subject we refer to the treatises of~\cite{Oja} and~\cite{DharChau}.

In the present paper, following an idea of the late Dietrich Morgenstern, we again start with the
marginal median as in the work of~\cite{Gr1996}, but then use a variant of the equipartition property 
instead of symmetrization. In dimension $d=2$, for example, a coordinate system
specifies four quarter spaces via the pairwise intersections of the half-spaces
underlying the marginal median, and a quarter median of $P$ may be defined
as any vector $\theta\in\bR^2$ such that for some coordinate system centered at $\theta$ 
all four of the respective quarter spaces have probability at least $1/4$. 
Figure~\ref{fig:staedte} shows a simple
example that Professor Morgenstern would have liked:  The data 
are the locations of the cities in Germany with a population of at 
least 100,000 (in 2017). We have $n=76$ such cities in total, so we can aim for a partition
with 19 cities in each open quarter space. One resulting quarter median is nearby the
exit Freudenberg of the A45 (Sauerlandlinie), and with the North-South axis tilted by about 
39 degrees counterclockwise; see the lines in Figure~\ref{fig:staedte}. The data (longitude/latitude)
for drawing the German border are extracted from the
file $\tt gadm36_{-}DEU_{-}0_{-}sp.rds$ available at {\tt https:\slash\slash gadm.org/data.html}.

Investigation of this concept leads to a number of interesting questions, starting with 
the obvious ones concerning existence and uniqueness. Next, what are the equivariance properties 
of the quarter median? Further, in contrast to the spatial median a quarter
median comes with a system of associated hyperplanes, which may be regarded as a part
of the estimator. In the two-dimensional case  these are just two orthogonal lines and can 
be parametrized by a counterclockwise rotation of the horizontal axis, thus providing a 
connection to a real interval.  
In higher dimensions, however, specification of a topological notion of
nearness of systems of coordinate hyperplanes is less straightforward, but it is needed
in connection with consistency of the estimators. For asymptotic normality even a local 
linearization is required. Finally, there is the algorithmic issue: How do we find a 
quarter median for a given data set?  

Section~\ref{sec:results} contains our main results,  on existence, equivariance,
and large sample behavior. Elliptical
distributions are treated in some detail and we compare the quarter median with other 
estimators for several specific families. 
We also consider algorithmic aspects together with the problem of measurable
selection from the respective solution sets. 
Proofs are collected in Section~\ref{sec:proofs}.

\section{Results}\label{sec:results} 

We first assemble some formal definitions and then give our results in separate subsections. 

\subsection{Preliminaries and basic definitions}\label{subsec:prelim}
We regard vectors $x\in\bR^d$ as column vectors and write $\langle x, y\rangle$ for 
the inner product of $x,y\in\bR^d$. The transpose of the vector $x$ is $x\t$,
so that $\langle x, y\rangle=x\t y$. Similarly, $A\t$ denotes the transpose of a matrix $A$.
We write $\|\cdot \|$ for the Euclidean norm on $\bR^d$,
$\Sd$ is the unit sphere in $\bR^d$, $\OO(d)$ is the group of orthogonal $d\times d$-matrices,
and $\SO(d)$ denotes the group of rotations, consisting of those elements of $\OO(d)$ that
have determinant~1. 

If $(\Omega,\cA,P)$ is a probability space, if $(\Omega',\cA')$ is a measurable space,
and if the mapping $S:\Omega\to\Omega'$ is $(\cA,\cA')$-measurable then we write $P^S$ for 
the push-forward of $P$ under~$S$, i.e.\ $P^S(A')=P(S^{-1}(A'))$ for 
all $A'\in\cA'$. If $(\Omega,\cA)=(\bR^d,\cB^d)$ and $b\in \Sd$ or 
$U\in \OO(d)$ then $P^{b}$ and $P^U$ respectively refer to the push-forwards 
of $P$ under the mappings $x\mapsto \langle b, x\rangle$ and $x\mapsto Ux$.

For $a\in\bR^d$ and $b\in \Sd$ we formally define 
the associated half-spaces by
\begin{equation*}
	 H_-(a;b) := \bigl\{x\in\bR^d:\, \langle b,x-a\rangle \le 0\bigr\}, \ 
	H_+(a;b):= \bigl\{x\in\bR^d:\, \langle b,x-a\rangle \ge 0\bigr\}.
\end{equation*}  
Let $e_1,\ldots,e_d$ be the vectors of the canonical basis of $\bR^d$. Then a
marginal median for a probability measure $P$
on $(\bR^d,\cB^d)$ is any vector $\theta\in\bR^d$ with the property that
\begin{equation}\label{eq:comeddef}
     P\bigl(H_+(\theta;e_i)\bigr)\ge 1/2,\  
     P\bigl(H_-(\theta;e_i)\bigr)\ge 1/2 , \quad 1\le i\le d.
\end{equation}
For two orthogonal unit vectors $b, b'\in \Sd$ and $a\in\bR^d$ the 
quarter spaces are the intersections of these half-spaces,
\begin{equation*}
	V_{\pm\pm}(a;b,b') := H_\pm(a;b)\cap H_\pm(a;b'),
\end{equation*}
with the obvious notational convention for the four possible combinations of plus and
minus signs. We then define a \emph{quarter median} of a probability measure $P$
on $(\bR^d,\cB^d)$ as any vector $\theta\in\bR^d$ with the property that, for some
orthonormal basis $B=\{b_1,\ldots,b_d\}$ of $\bR^d$,
\begin{equation}\label{eq:qmeddef1}
  P\bigl(H_+(\theta;b_i)\bigr)\ge 1/2,\  P\bigl(H_-(\theta;b_i)\bigr)\ge 1/2 , \quad 1\le i\le d,
\end{equation}  
and  
\begin{equation}\label{eq:qmeddef2}
                  P\bigl(V_{\pm\pm}(\theta;b_i,b_j)\bigr)\ge 1/4, \quad 1\le i<j\le d,
\end{equation}
for all combinations of plus and minus signs. 
If~\eqref{eq:qmeddef1} and~\eqref{eq:qmeddef2} hold then we call the pair 
$(\theta,U)\in\bR^d\times \OO(d)$, where $U$ has rows $b_1\t,\ldots,b_d\t$,
a solution of the \emph{quarter median problem} for the probability measure $P$. 
Again, the data versions of these notions refer to the empirical distribution $P_n$.  

If $P$ is absolutely continuous then hyperplanes have zero probability and
\eqref{eq:qmeddef1} follows from \eqref{eq:qmeddef2}. In general, this is not the case 
as the following simple example shows: Let $d=2$ and
\begin{equation*}
          P=\tfrac{1}{3}\bigl(\delta_{(0,0)\t}+\delta_{(1/2,1)\t}
              +\delta_{(1,1/2)\t}\bigr).
\end{equation*}
Then \eqref{eq:qmeddef2} holds with $b_1=e_1$, $b_2=e_2$ and $\theta = (1,1/2)\t$, 
but \eqref{eq:qmeddef1} is not true. Further, with this choice of 
$b_1,b_2$ it is easily verified that the unique marginal median 
$\theta=(1/2,1/2)\t$ is a quarter median.
This example also shows that the quarter median may not be unique.
Further, 
if hyperplanes have probability zero, then the inequalities in \eqref{eq:qmeddef1} and 
\eqref{eq:qmeddef2} can be replaced by equalities, and with $(\theta,U)$ as solution 
of \eqref{eq:qmeddef2} we then speak of a solution of the \emph{equipartition problem} 
for $P$.

\subsection{Existence and uniqueness}\label{subsec:exist}
The existence of a median and thus of a marginal median follow immediately from the 
monotonicity of distribution functions. For the quarter median counting the number of 
variables and the number of constraints may give a first impression. For simplicity, 
we temporarily assume that hyperplanes have probability zero, so that we have an 
equipartition problem. If $d=2$, there 
are then three unknowns, two for the location parameter and one for the angle of rotation,
see also Remark~\ref{rem:uniqueness}\,(b) below,
and one of the four constraints is redundant, so that the number of unknowns is the same 
as the number of equations. For $d>2$ we first note that the conditions 
in~\eqref{eq:qmeddef1} and~\eqref{eq:qmeddef2} are not independent. Indeed, 
for a given basis $\{b_1,\ldots,b_d\}$ of $\bR^d$ it is enough to have probability  
$1/2$ for $H_+(\theta,b_i)$, $i=1,\ldots,d$, and probability $1/4$ for $V_{++}(\theta;b_i,b_j)$
for $1\le i < j\le d$, so that $d + d(d-1)/2$ conditions remain. On the other hand, we have 
$d(d-1)/2$ parameters for the (special) orthogonal group, and there are $d$ components of
the location vector, hence the number of unknowns is again the same as the number of independent
constraints. This, of course, is just a heuristic argument, as the respective conditions are nonlinear.

Using a general tool from algebraic topology, the homotopy invariance of the Brouwer
degree of continuous mappings,~\cite{Mak2007} solved the equipartition
problem under the additional condition that $P$ is absolutely continuous. 
A smoothing argument leads to the extension needed here. 

\begin{theorem}\label{thm:main2}
A solution $(\theta,U)$ of the quarter median problem exists for every probability measure
$P$ on $(\bR^d,\cB^d)$.
\end{theorem}

This result raises an obvious question: 
Is there a stronger equipartition in dimension $d>2$, such as an 
`octomedian' if $d=3$, with each of the octants receiving the same probability?
This is certainly possible for specific distributions, such as the multivariate
standard normal. However,  if $d=3$ then we have three
unknowns for the location vector, and three for the rotation, which can be
parametrized by the Euler angles, and the condition that the octants all have
probability $1/8$ would similarly lead to seven constraints, which is one too many.
Obviously, the difference would be even greater for dimension $d>3$ if
all $2^d$ intersections of coordinate half-spaces were to have the same probability.
According to~\cite[p.554]{Mak2007} this dimension counting argument already 
implies that a solution does not exist for a generic distribution. 

The notion of uniqueness requires some attention. We say that a distribution $P$
has a unique quarter median if $\theta_1=\theta_2$ for any two solutions 
$(\theta_1,U_1)$ and $(\theta_2,U_2)$ of the quarter median problem. 
Obviously, the full pair
$(\theta,U)$ will never be unique: Any linear hyperplane, i.e.~a subspace $H$ of $\bR^d$
of dimension $d-1$, may be specified as the orthogonal complement of the one-dimensional 
space (line) generated by a unit vector $b\in\Sd$ via $H(b):=\{x\in\bR^d:\, b\t x=0\}$, 
but obviously $-b$ would lead to the same hyperplane. 
This simple observation already
implies that we may restrict $U$ to be an element of $\SO(d)$ when searching for  
a quarter median. 
Further, the conditions~\eqref{eq:qmeddef1} and~\eqref{eq:qmeddef2} are invariant under
permutations of the rows of $U$. 
Conversely, given a set of coordinate hyperplanes, a corresponding basis 
is a set with elements that are unique up to a factor $-1$. Putting this together
we write $U\sim V$ for $U,V\in\OO(d)$ if $U$ and $V$ lead to the same set of linear
hyperplanes; we then say that the solution of the quarter median problem is unique
if for any two solutions $(\theta_1,U_1)$ and $(\theta_2,U_2)$ 
we have $\theta_1=\theta_2$ and~$U_1\sim U_2$.

The following formal approach will turn out to be useful. 
Let $G_1$ be the subgroup of $\OO(d)$ that consists of the diagonal matrices with
diagonal entries from $\{-1,+1\}$, i.e. the subgroup that represents the compositions
of reflections at the coordinate axes, and let $G_2$ be the subgroup of 
permutation matrices, which corresponds to coordinate permutations.  
We write $G$ for the subgroup generated by $G_1$ and $G_2$. This is the system
of all $d\times d$-matrices $A$ 
with entries $a_{ij}\in\{-1,0,1\}$ and such that for some permutation $\pi$ of $\{1,\ldots,d\}$ 
we have  $a_{ij}\not=0$ if and only if $j=\pi(i)$. The above transition from equality to 
equivalence can be then regarded as 
a transition from the group $\OO(d)$ to the factor group
\begin{equation}\label{eq:OOtoH}
             \bH(d) := \OO(d)/G,
\end{equation} 
with the equivalence classes $[U] :=\{V\in\OO(d):\, V\sim U\}$, $U\in\OO(d)$, as elements. 
If we think of a solution of the quarter median problem as a pair 
$(\theta,H)\in\bR^d\times \bH(d)$, then uniqueness means that two solutions are equal.

\begin{remark}\label{rem:uniqueness}
(a) The $\theta$-part of a solution $(\theta,H)$ for the quarter median problem may
be unique while the $H$-part is not. For example, if $P$ has a density that depends 
on $x\in\bR^d$ only through $\|x\|$ then $0$ is the unique quarter median, but all 
$H\in\bH(d)$ lead to a solution $(0,H)$. 

(b) In dimension $d=2$ there is an (almost) canonical choice of an element $U\in\OO(d)$ 
from an equivalence class $H\in\bH(d)$ corresponding to a set of hyperplanes, here 
just two orthogonal lines $L_1,L_2$ meeting at the point $0\in\bR^2$. With a suitable 
half-open interval $I$ of length $\pi/2$, such as $I=[-\pi/4,\pi/4)$, there is exactly one 
$\alpha\in I$ such that the lines are represented by    
\begin{equation}\label{eq:parSO2}
	U_\alpha:=\begin{pmatrix} \cos(\alpha)&-\sin(\alpha)\\ 
		\sin(\alpha) & \cos(\alpha)\end{pmatrix},
\end{equation}
which is the counterclockwise rotation of the canonical coordinate lines by the angle~$\alpha$,
in the sense that $L_i=\{x\in\bR^2:\, b_i\t x=0\}$, where $b_1\t,b_2\t$ are the rows of $U$. 
The corresponding topology would then be that of the quotient space $\bR/I$.

(c) While we work with the space $\bR^d$ of column vectors throughout the notions of
quarter median and quarter median problem make sense for an arbitrary finite-dimensional
Euclidean space $(E,\langle\cdot,\cdot \rangle)$. 
If $(\theta,U)$ solves the quarter median problem for a distribution $P$
on $(\bR^d,\cB^d)$, then $\{b_2,\ldots,b_d\}$ is an orthonormal basis for the linear space 
$E=H(b_1)$, which has dimension $d-1$. Here $b_i\t$ denotes the $i$th row of $U$, $i=1,\ldots,d$;
equivalently, $b_i$ is the $i$th column of $U\t$. If $d\ge 3,$ this can be used to obtain a
projectivity or consistency property that relates 
$(\theta,U)$  to a solution $(\theta_E,U_E)$ of the quarter median problem
for the push-forward of $P$ under the orthogonal projection $\pi_E:\bR^d\rightarrow E$.
In fact, let $\Sd^E=\{y\in E:\|y\|=1\}$ be the unit sphere in $E.$ For $c\in E$ and orthogonal
unit vectors $d,d'\in \Sd^E$ the associated
half-spaces and quarter spaces are given by
\begin{equation*}
  H_-^E(c;d) := \bigl\{y\in E:\, \langle d,y-c\rangle \le 0\bigr\}, \ 
  H_+^E(c;d) := \bigl\{y\in E:\, \langle d,y-c\rangle \ge 0\bigr\},
\end{equation*}
and
\begin{equation*}  
  V_{\pm\pm}^E(c;d,d') := H_\pm^E(c;d)\cap H_\pm^E(c;d').
\end{equation*}
Using that a vector $x\in\bR^d$ can uniquely be written as 
$x=\sum_{i=1}^d\langle x, b_i \rangle\,b_i$ we see that the orthogonal projection
of $x$ in $E$ is $x_E:=\pi_E(x)=\sum_{i=2}^d\langle x, b_i \rangle\,b_i$. 
It follows from
$\pi_E^{-1}\left (H_{\pm}^E(\theta_E;b_i)\right ) = H_{\pm}(\theta;b_i)$ for $2\le i\le d$,
and
$\pi_E^{-1}\left (V_{\pm\pm}^E(\theta_E;b_i,b_j)\right ) = V_{\pm\pm}(\theta;b_i,b_j)$ for $2\le i<j\le d$, 
that
$P^{\pi_E}\left (H_{\pm}^E(\theta_E;b_i)\right )=P\left (H_{\pm}(\theta;b_i)\right )$
for $2\le i\le d$, and that
$P^{\pi_E}\left (V_{\pm\pm}^E(\theta_E;b_i,b_j)\right )=P\left (V_{\pm\pm}(\theta;b_i,b_j)\right )$
for $2\le i<j\le d$. 
\end{remark}

\subsection{Equivariance}\label{subsec:equi}
In the context of equivariance properties of location estimators for $d$-dimensional
data we start with a family (often a group) $\cS$ of measurable transformations 
$S:\bR^d\to\bR^d$ and we regard the location parameter $\theta\in\bR^d$ as a function 
$\theta=T(P)$ of $P$.  Recall that $P^S$ is 
the push-forward of $P$ under the transformation $S$. Then $T$ is
said to be \emph{$\cS$-equivariant}
if 
\begin{equation*}
   T(P^S)=(S\circ T)(P)\quad\text{for all }S\in\cS.
\end{equation*} 
In cases where the parameter
is not unique this might better be expressed as $S(\theta)$ being a parameter of this type
for $P^S$ if $\theta$ is for $P$. 

For example, with $\cS$ being one of the groups $G_1$ and $G_2$ that were used in connection 
with~\eqref{eq:OOtoH}, it is easy to check that
\begin{equation}\label{eq:Ginv}
    S^{-1} \bigl(\mMed_m(P^S)\bigr) = \mMed_m(P) \quad\text{for all } S\in\cS,
\end{equation} 
i.e.~the specific marginal median introduced in \eqref{eq:midpointd} is equivariant with respect to reflections and coordinate
permutations. It is easily seen that we may even take $\cS=G$. This specific marginal median 
is also equivariant with respect to separate rescalings, where the corresponding
group is the set of all regular diagonal matrices. 

The general idea is of course that statistical inference should respect the structure 
of the data space. Equivariance of a location estimator with respect to shifts and
orthogonal linear transformations means that the estimator interacts in this way
with the full isometry group of Euclidean space.

At various stages the following function will be important,
\begin{equation}\label{eq:defPhi}
    \Psi(U,P)\, :=\, U\t\mMed_m(P^U),
\end{equation}
with $U\in\OO(d)$ and $P$ a probability measure on $(\bR^d,\cB^d)$; $\Psi(U,P)$ is the
specific marginal median in the coordinate system associated with $U$, transformed back to
canonical coordinates. By definition, a quarter median is always a marginal median in
a specific coordinate system. Below, marginal medians in 
\begin{equation}\label{eq:MMMP}
  \cM(P)\; :=\; \bigtimes_{i=1}^d
                \Bigl\{\Med_-(P^{\pi_i}),\Med_m(P^{\pi_i}),\Med_+(P^{\pi_i})\Bigr\}
\end{equation}
will be of particular interest. 

\begin{proposition}\label{prop:qumedequi}
Let $P$ be a probability distribution on $(\bR^d,\cB^d)$, and let $\cS$ be the set of
Euclidean motions.

\vspace{.2mm}
\noindent
\emph{(a)} If $\theta$ is a quarter median for the distribution $P$ then, for all $S\in\cS$,
$S(\theta)$ is a quarter median for $P^S$.

\vspace{.5mm}
\noindent
\emph{(b)} If $U,V\in\OO(d)$ are such that $U\sim V$ then $\Psi(U,P)=\Psi(V,P)$. 

\vspace{.5mm}
\noindent
\emph{(c)} If $(\theta,U)$ solves the quarter median problem for $P$, then
there is a marginal median $\eta \in \cM(P^U)$ such that
$(\Tilde \theta,U)$ with $\Tilde\theta=U\t \eta$ is also a solution.    

\vspace{.5mm}
\noindent
\emph{(d)} If $(\theta,U)$ solves the quarter median problem for $P$ and $S\in\cS$,
$S(x)=Ax+b$, then $(A\theta + b, U\!A\t)$ solves the quarter median problem for $P^S$.
Moreover, if $P^b$ has a unique median for all $b \in \Sd$,
\begin{equation}\label{eq:equiselect}
   \Psi\bigl(U\!A\t,P^S\bigr)\, = \, S\bigl(\Psi(U,P)\bigr).
 \end{equation}

\vspace{.5mm}
\noindent
\emph{(e)} Suppose that $P^b$ has a unique median for all $b \in \Sd$. Then, the function
$\Psi(\,\cdot\,,P):\OO(d)\to \bR^d$ is continuous.
\end{proposition}
 
Part (a) is the equivariance on the set level. Part (b) shows that $\Psi(\cdot,P)$
may be regarded as a function on $\bH(d)$. Part (c) is of importance when developing
an algorithm for searching a solution of empirical quarter median problems: 
For each $U\in\OO(d)$ the set of potential $\theta$-values in~\eqref{eq:MMMP} 
consists of three
vectors only, separately for each coordinate. 
It is tempting to think that this could be reduced further to the 
respective midpoint, but the uniform distribution on the six points 
$(-2,-2)\t,(-1,3)\t,(1,-1)\t,(2,2)\t,(3,4)\t,(4,0)\t$ in the plane provides 
a counterexample: The pair $\bigl((2,2)\t,U\bigr)$ with $U=\diag(1,1)$ solves the quarter
median problem for this distribution, but the midpoint of the 
associated rectangle $[1,2]\times[0,2]$ of multivariate medians 
would lead to $\bigl((1.5,1)\t,U\bigr)$, which is not a solution.
Part (d) notices the Euclidean motion equivariance of the quarter median
for distributions $P$ with the property that $P^b$ has a unique median for all $b\in\Sd$,
a condition that will be used repeatedly below.

Similar to the spatial median and the orthomedian introduced by~\cite{Gr1996}, the quarter 
median is not affine equivariant. In particular, it is not equivariant with respect to separate
rescalings of the coordinates, in contrast to the special marginal median, but
for distributions $P$ with the property that $P^b$ has a unique median for all $b\in\Sd$,
equivariance does hold (and is used in the proofs below) for simultaneous rescalings, where the group $\cS$ 
consists of all positive scalar multiples of the identity matrix. 

\begin{figure}
\vspace{-.2cm}
\setlength{\abovecaptionskip}{-0.2cm}
\hskip1cm\includegraphics[scale=.5]{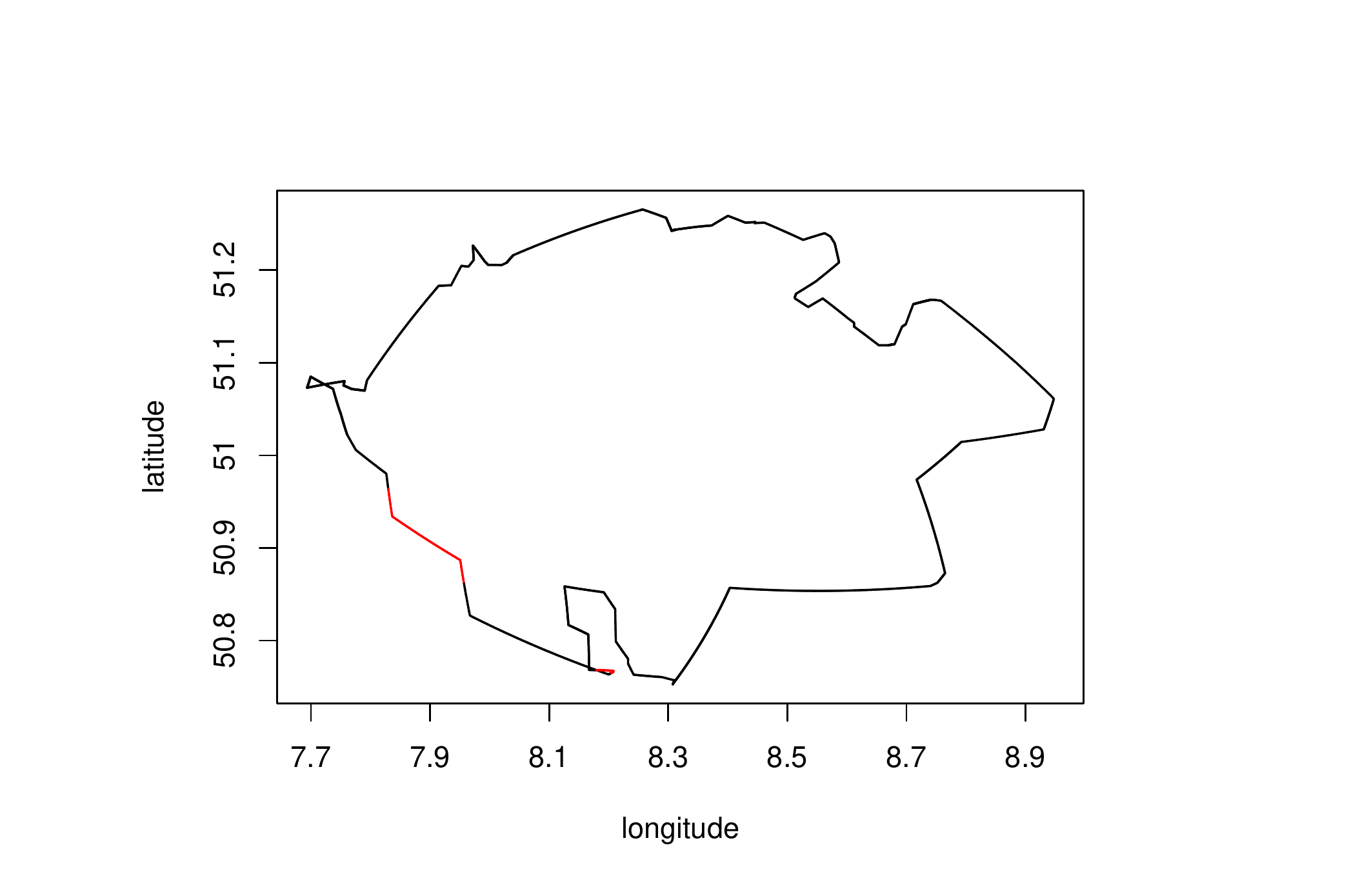}
\caption{The median curve  for the city data}
\label{fig:staedteMMed}
\end{figure}

Any quarter median is a marginal median in a suitably chosen coordinate
system, which we may take to be generated by an orthogonal matrix with determinant~1. 
\cite{Gr1996} introduced the orthomedian of $P$ as the integral of the function
$\Psi(\cdot,P)$ over $\SO(d)$ with respect to the Haar measure with
total mass~1. Thus, the orthomedian may be seen 
as the expected value of this function with $U$ chosen uniformly at random, 
and its orthogonal equivariance is then an immediate consequence
of the fact that the Haar measure is invariant under the group operations. In contrast,
the quarter median picks a marginal median from the range of the set-valued function
\begin{equation*}
 U\mapsto \Bigl\{U\t \eta: \, \eta\in \bigl[\mMed_-(P^U),\mMed_+(P^U)\bigr]\Bigr\}
\end{equation*}
via the partition constraint~\eqref{eq:qmeddef2}. Figure~\ref{fig:staedteMMed} shows the range of the function
$U\mapsto\Psi(U,P_n)$ for the empirical distribution $P_n$ associated with the
city data in Figure~\ref{fig:staedte}, where $n=76$. The parts colored red correspond to possible quarter medians.

\subsection{Asymptotics}\label{subsec:asymp}
Throughout this section we assume that $P$ is a probability 
measure on $(\bR^d,\cB^d)$ with the property that
\begin{equation}\label{as:cond1}
  P^a~\text{has a unique median for all}~a \in S_d,\hskip3.3cm
\end{equation}
and that $X_1,X_2,\ldots$ are independent copies of a random vector $X$ with distribution $P$.
We further assume that, for each $n\in\bN$, the pair $(\theta_n,U_n)\in \bR^d\times \OO(d)$ 
with  $\theta_n=\theta_n(X_1,\ldots,X_n)$, $U_n=U_n(X_1,\ldots,X_n)$, is a solution of the
quarter median problem for the (random) empirical distribution $P_n$ associated with
$X_1,\ldots,X_n$. We also assume that $\theta_n$ and $U_n$ are measurable with respect 
to the respective Borel $\sigma$-fields on $(\bR^d)^n$. A discussion of 
the associated selection problem and an extension of the results to more general
random elements is given in Section~\ref{subsec:select}. 

We deal with consistency first and then consider asymptotic normality.
Recall that uniqueness and convergence refer not to the $U$-matrices themselves,
but to the associated equivalence classes $[U]\in\bH(d)$,
hence we have to specify what convergence 
of the $U$-part means. On $\OO(d)$ we use the topology induced by $\bR^{d\times d}$
where a sequence of matrices converges if all entries converge individually. It is
well known that this makes $\OO(d)$ a compact and separable Hausdorff space. 
We use the quotient topology on $\bH(d)$, which is  
the finest topology with the property that the mapping $\SO(d)\to\bH(d)$, 
$U\mapsto [U]$, is continuous. As we have
$\bH(d)=\OO(d)/G$ with some finite group $G\subset \OO(d)$,
the space $\bH(d)$ with this topology is again a compact and separable Hausdorff space. 

\begin{theorem}\label{thm:consistency} Suppose that  \eqref{as:cond1} holds.

\noindent
\emph{(a)} If the quarter median for $P$ is unique and given by
$\theta\in\bR^d$, then  $\theta_n$ converges almost surely to $\theta$ 
as $n\to\infty$. 

\vspace{.3mm}

\noindent
\emph{(b)} If  the solution of the quarter median problem for $P$ is unique and 
given by $(\theta,H)\in \Sd\times\bH(d)$, then $(\theta_n,[U_n])$ 
converges almost surely to $(\theta,H)$ as $n\to\infty$. 
\end{theorem} 

In connection with distributional asymptotics we require~\eqref{as:cond1} 
and also that $P$ is absolutely continuous. Let $M_a=M(a)$ be
the (unique) median of $P^a$, $a\in \Sd$.   
As in the paper of~\cite{Gr1996} we further assume that $P^a$ has a density $f_a$ such that
\begin{align}
  &(a,t)\rightarrow f_a(t)\text{ is continuous at each}~(a,t)
     \in        D:=\big \{\big (a,M_a\big ):a\in \Sd\big \},\label{as:cond2}\\
  &f_a(M_a)>0~\text{for each}~a\in\Sd, \label{as:cond3}\
\end{align}
and, in addition, that 
\begin{align}
  &\text{there is a point}~\theta\in\bR^d~\text{such that }M_a=a\t\theta
   \text{ for each}~a\in\Sd.\label{as:cond4}
\end{align}
For example, \eqref{as:cond4} is satisfied if $P$ is symmetric about $\theta$
in the sense that $(X_1-\theta)$ and $-(X_1-\theta)$ have the same distribution. 
Of course, if \eqref{as:cond1} and \eqref{as:cond4} hold then the quarter median 
for $P$ is unique and given by $\theta$.

For $U\in\OO(d)$ with rows $b_1\t,\ldots,b_d\t\,$ let 
\begin{equation}\label{eq:defDelta}
  \Delta(P,U) \, :=\, 
     \diag\biggl(\frac{1}{f_{b_1}(M_{b_1})^2},\ldots, \frac{1}{f_{b_d}(M_{b_d})^2}\biggr).
\end{equation} 
Of course, $f_a$ and $M_a$ also depend on $P$. Finally, we write $X_n\dconv Q$ for
convergence in distribution of a sequence $(X_n)_{n\in\bN}$ of random variables to a 
random variable with distribution $Q$, and $N_d(\mu,\Sigma)$ for the
$d$-dimensional normal distribution with mean vector $\mu$ and covariance matrix~$\Sigma$.

\begin{theorem}\label{thm:CLTtheta}
  Suppose that the conditions \eqref{as:cond1}--\eqref{as:cond4} hold and that the solution of the
  quarter median problem for the absolutely continuous probability measure $P$ is unique and given
  by $(\theta,H)$. Then, in the above setup, and with $U\in H$,
  \begin{equation}\label{eq:CLTtheta}
    \sqrt{n}\bigl(\theta_n-\theta) \, \dconv \, N_d\bigl(0,\Sigma(P,U)\bigr)
                              \quad\text{as } n\to\infty, 
  \end{equation}
  where $\,\Sigma(P,U):= U\t\Delta(P,U)U$ and $\Delta(P,U)$ is given by~\eqref{eq:defDelta}.
\end{theorem}

As part of the proof we will show that $\Sigma(P,U)$ does not depend on the choice of the 
element $U$ of $H$.

\subsection{Elliptical distributions}

We now consider a special family of distributions in more detail. 
Let $h:\bR_+\to\bR_+$ be a function
with the property
\begin{equation}\label{eq:defcd}
     0 \, <\, c_d(h):= \int_{\bR^d} h(\|x\|^2)\, dx \, <\, \infty.
\end{equation}
Then, for each $\mu\in\bR^d$ and each positive definite matrix $\Sigma\in\bR^{d\times d}$,
\begin{equation*}
    f(x) = f(x;\mu,\Sigma)\; =\; \frac{1}{c_d(h) (\det \Sigma)^{1/2}} 
                          \,h\bigl((x-\mu)\t\Sigma^{-1}(x-\mu)\bigr), \quad x\in\bR^d,
\end{equation*}
is the density of a probability measure $P$ on $(\bR^d,\cB^d)$, the 
\emph{elliptical distribution} 
with \emph{location vector} $\mu$ and \emph{dispersion matrix} $\Sigma$.
We abbreviate this to $P=\Ell_d(\mu,\Sigma;h)$. If $\mu=0$ and $\Sigma=I_d$ we write
$\Sym_d(h)$ for $\Ell_d(0,I_d;h)$ and speak of a \emph{spherically symmetric} 
distribution with the
density $h\big (\| x \|^2\big )/c_d(h),\,x\in \bR^d$. A transformation to
spherical coordinates leads to
\begin{equation}\label{eq:defcdalt}
c_d(h)\, =\, \frac{2\pi^{d/2}}{\Gamma(d/2)} \int_0^\infty r^{d-1} h(r^2)\, dr  =\, \frac{\pi^{d/2}}{\Gamma(d/2)} \int_0^\infty r^{\frac{d}{2}-1} h(r)\, dr.
\end{equation}
The function $h$ is said to be the \emph{density generator} of the elliptical distributions $\Ell_d(\mu,\Sigma;h).$
For $h(t)=\exp(-t/2)$ the associated constant is $c_d(h)= (2\pi)^{d/2}$, and we obtain the 
multivariate normal distributions $N_d(\mu,\Sigma)$.

For later purposes we list some properties of elliptical distributions. Clearly, 
if $P=\Ell_d(\mu,\Sigma;h)$ 
and if $T:\bR^d\to\bR^d$, $T(x) = Ax + b$ with $A\in\bR^{d\times d}$
regular, $b\in\bR^d$, then $P^T=\Ell_d(\mu+b, A \Sigma A\t; h)$. Further, for $d\ge 2$ the 
univariate marginal densities of $P=\Sym_d(h)$ coincide and are given by
\begin{equation*}
       g_1(y)=g(y^2)/c_d(h),\ y\in \bR,
\end{equation*}
where
\begin{align*}
  g(y)\ &=\ \int_{\bR^{d-1}} h(y + x_2^2+ \cdots + x_d^2) \, dx_2\cdots dx_d\\
       &=\ \frac{\pi^{\frac{d-1}{2}}}{\Gamma\big (\frac{d-1}{2}\big )}\int_{y}^\infty h(t)\big (t-y \big )^{\frac{d-1}{2}-1}\,dt,\quad y\ge 0.
\end{align*}
The function $g$ is nonincreasing and continuous on $(0,\infty)$. Hence the univariate
marginal distributions are unimodal with mode 0 and their density is continuous on
$\bR\setminus \{0\}$; see e.g.\ \cite[p. 37]{Fang}. 

In what follows we deal with generators $h$ that satisfy the condition
\begin{equation}\label{cond:stetig}
  \sup_{0\le t\le \epsilon} h(t)<\infty\hskip3mm\text{for some finite interval}~[0,\epsilon].
\end{equation}
Then $g$ and $g_1$ are bounded and continuous. Additionally, with
\begin{equation}\label{eq:defcdgen}
c_k(h)\, 
  :=\, \frac{2\pi^{k/2}}{\Gamma(k/2)} \int_0^\infty r^{k-1} h(r^2)\, dr \, 
 =\, \frac{\pi^{k/2}}{\Gamma(k/2)} \int_0^\infty r^{\frac{k}{2}-1} h(r)\, dr
\end{equation}
for $k=1,\ldots,d-1$, we have $0<c_k(h)<\infty$; especially, $g(0)=c_{d-1}(h).$ 
Finally, the transformation law mentioned above can be generalized to $k\times d$-matrices with
rank $k<d$ at the cost
of changing the function $h$. In particular, the first marginal distribution associated with
$\Ell_d(0,\Sigma;h)$ is $\Ell_1(0, \sigma_{11}; g)$ if $\Sigma=(\sigma_{ij})_{i,j=1}^d$,
and has a density that is bounded, continuous and strictly positive at 0.
More generally, if $P=\Ell_d(\mu,\Sigma;h)$ then the density
$f_a$ of $P^a$ for $a\in\Sd$ is given by
\begin{align*}
  f_a(t)\,=\,\frac{1}{c_d(h) (a\t \Sigma a)^{1/2}}\,g\left (\frac{(t-a\t \mu)^2}{a\t \Sigma a}\right ), \quad t\in\bR,
\end{align*}
$M_a=a\t \mu$ is the unique median of $P^a$, $f_a(M_a)$ is positive for all $a\in\Sd,$ and the
function $(a,t)\rightarrow f_a(t)$ is continuous. Hence the assumptions
\eqref{as:cond1} - \eqref{as:cond4} of Theorem~\ref{thm:CLTtheta} are satisfied.
 
As $\Sigma$ is positive definite and symmetric
the principal axes theorem applies and gives the representation
\begin{equation}\label{eq:priax}
   \Sigma \,=\, U\t\, \diag(\lambda_1,\ldots,\lambda_d)\, U
\end{equation}  
with some $U\in\OO(d)$ and $\lambda_1\ge\lambda_2\ge\cdots\ge\lambda_d>0$, where the
$\lambda_i$'s are the eigenvalues of $\Sigma$. The distribution $P=\Ell_d(\mu,\Sigma;h)$
is said to be \emph{strictly elliptical} if no two of these eigenvalues coincide.
The next theorem deals with the uniqueness of (the solution of)
the quarter median (problem) for
elliptical distributions. Note that condition \eqref{cond:stetig} is 
not imposed there.

\begin{theorem}\label{thm:existell}
If $P=\Ell_d(\mu,\Sigma;h)$,
then the quarter median of $P$ is unique and given by~$\mu$.
Moreover, if $P$ is strictly elliptical, then the solution of the quarter median
problem is unique and given by $(\mu,[U])$ with $U$ as in~\eqref{eq:priax}.
\end{theorem}

For samples from strictly elliptical distributions with density generators satisfying
\eqref{cond:stetig} we obtain asymptotic normality for a
sequence of associated quarter medians $\theta_n=\theta_n(X_1,\ldots,X_n)$.
Part (a) of the following result is essentially
a corollary to Theorem~\ref{thm:CLTtheta}. In the symmetric case we 
only have weak uniqueness, hence we need to consider this situation
separately. Intermediate cases, where only
some of the eigenvalues coincide, can be treated similarly, but require 
a certain amount of bookkeeping.

\begin{theorem}\label{thm:CLTell}
Let $X_n$, $n\in\bN$, be independent random vectors with distribution $P=\Ell_d(\mu,\Sigma; h)$, 
with $h$ satisfying \eqref{cond:stetig}. Let $c_{d-1}(h),c_d(h)$ be as in~\eqref{eq:defcd}
and let 
\begin{equation*}
          \sigma_{\QMed}^2 :=  c_{d}(h)^2/(2c_{d-	1}(h))^2.
\end{equation*}

\emph{(a)} Suppose that $(\theta_n,U_n)_{n\in\bN}$ is a sequence of random variables 
with values in $\bR^d\times \OO(d)$
such that, for all $n\in\bN$,  $(\theta_n,U_n)$ solves the quarter median
problem for~$P_n$.
Then, if $P$ is strictly elliptical,
\begin{equation}\label{eq:CLT}
  \sqrt{n}(\theta_n-\mu)  \; \dconv \; N_d(0,\sigma_{\QMed}^2\,\Sigma)\quad\text{as } n\to\infty.
\end{equation}

\emph{(b)} Suppose that $\Sigma=\lambda I_d$ for some $\lambda>0$. Let $(U_n)_{n\in\bN}$ be an 
arbitrary sequence of elements of $\OO(d)$ and suppose that, for all $n\in\bN$,
$\theta_n$ is a marginal median of~$P_n^{U_n}$. Then
\begin{equation}\label{eq:CLTsymm}
  \sqrt{n}(\theta_n-\mu)  \; \dconv \; N_d(0,\sigma_{\QMed}^2\,\Sigma)\quad\text{as } n\to\infty.  
\end{equation}
\end{theorem}

Remarkably, the covariance matrix of the limiting normal distribution always is
a fixed multiple of the dispersion matrix $\Sigma$. This is in contrast to other
familiar estimators of $\mu$ such as the empirical spatial median $\SMed(X_1,\dots,X_1),$
for example. In fact, with reference to \cite{Bai} where the asymptotics of the
spatial median are considered in an even more general context,~\cite{Somorcik}
states that in the $d$-dimensional case under some weak conditions 
\begin{equation}\label{smed:asym}
  \SMed(X_1,\dots,X_1) \; \dconv \; N_d(0,V),
\end{equation}
where the asymptotic covariance matrix is given by
\begin{equation}\label{smed:asympcov}
  V=D_1^{-1}D_2D_1^{-1}
\end{equation}
with
\begin{equation*}
  D_1=E\left (\frac{1}{\| X_1-\mu \|}\left (I_d - \frac{(X-\mu)(X-\mu)\t}{\|X-\mu\|^2}\right )\right )
\end{equation*}
and
\begin{equation*}
   D_2=E\left (\frac{(X-\mu)(X-\mu)\t}{\|X-\mu\|^2}\right ).
\end{equation*}
In the special (spherical) symmetric case $\Sigma=I_d$ the identity \eqref{smed:asympcov} simplifies to
\begin{equation*}
  V=\sigma_{\SMed}^2I_d,
\end{equation*}
where
\begin{equation}\label{smed:asymcov:spez}
  \sigma_{\SMed}^2=\frac{d}{(d-1)^2}\left (E\left (\frac{1}{\| X-\mu\|}\right )\right )^{-2}.
\end{equation}
It was shown by~\cite{Gr1996} that in this symmetric case the empirical orthomedian has the same
asymptotic covariance matrix as the spatial median. 
In what follows, we give a comparison of the covariances of the limit distributions
for some estimators of $\mu$, 
in particular for $\QMed(X_1,\dots,X_n)$ as the empirical quarter median, the
sample mean $\frac{1}{n}\sum_{j=1}^nX_j$, and the maximum likelihood estimator $\ML(X_1,\dots,X_n)$, 
for some special distribution families, where our main focus is on strictly elliptical distribution families. In all cases, the limit distribution of the standardized 
estimators is a centered multivariate normal, and by asymptotic covariance we mean the respective 
covariance matrix.

\begin{example}\label{ex:CLT} 
We work out the details for some specific distribution
families in the case of dimension $d=2$. 

\smallbreak 
(a) For normal distributions $P=N_2(\mu,\Sigma)$ Theorem~\ref{thm:CLTell} (a) leads to
\begin{equation}\label{eq:CLTnormal}
  \sqrt{n}(\QMed (X_1,\dots,X_n)-\mu)  \; \dconv \; N_2(0,\tfrac{\pi}{2}\Sigma).
\end{equation}
The maximum likelihood estimator of $\mu$ is the sample mean; its (asymptotic) covariance
is simply $\Sigma$ itself, hence the asymptotic efficiency of the quarter median 
is about $64\%$. Of course, the quarter median is more robust with respect to gross outliers
in the data. In the symmetric case $\Sigma=I_2$ both orthomedian and spatial median have
asymptotic covariance $\tfrac{4}{\pi}I_2$. 

For non-symmetric two-dimensional normal distributions with $\Sigma=\diag(1,\lambda)$,
$0<\lambda<1$,  \cite{BrownSMed} and~\cite{Gr1996} obtained expressions for the
asymptotic covariance matrices of the spatial median and the orthomedian. 
Both are of diagonal form, with the two diagonal entries
given by expressions involving infinite series and double integrals respectively. These
can be used numerically, see~\cite[Table 1]{BrownSMed} and~\cite[Table 1]{Gr1996}. For the 
second value, referring to the shorter axis, we have the simple explicit formula 
$\lambda\pi/2$ for the quarter median,
which is smaller than the value for the spatial median if  $\lambda<0.038$, and 
for the orthomedian if~$\lambda<0.1$.

\smallbreak
(b) Taking $h(t)=\exp(-\tfrac{1}{2}\sqrt{t}),\,t\ge 0,$ we obtain the bivariate doubly exponential 
distributions ${\rm MDE}_2(\mu,\Sigma)$, with densities given by
\begin{equation*}
  f(x;\mu,\Sigma)=\frac{1}{8\pi}\frac{1}{\sqrt{\det \Sigma}}
  \exp\left (-\tfrac{1}{2}\left ((x-\mu)\t \Sigma^{-1}(x-\mu)\right )^{\frac{1}{2}}\right ),~x\in\bR^2.
\end{equation*}
The asymptotic covariance matrix for the sample mean of independent two-dimensional
random vectors $X_1,\dots,X_n$ with the density $f(\cdot;\mu,\Sigma)$ is given by $\cov(X_1)=12\,\Sigma$;
see, e.g. \cite[Proposition 3.2 (iii)]{Gomez}.  Using
\begin{equation*}
   c_1(h)=2\int_0^\infty  e^{-\frac{1}{2}t}\,dt
              =4,\hskip5mm  c_2(h)=2\pi \int_0^\infty t e^{-\frac{1}{2} t}\,dt=8\pi,
\end{equation*}
we obtain $\sigma^2_{\QMed}=\pi^2$, so that $\pi^2\Sigma$ is the asymptotic covariance matrix for
the sample quarter median. For a given fixed symmetric positive
definite $\Sigma$ the Fisher information matrix can be calculated to be $\frac{1}{8}\Sigma^{-1}$
(see, e.g. \cite[formulas (3.2), (3.3)]{Mitchell}). So, by the general asymptotic theory of
maximum likelihood estimation $8\Sigma$ is the asymptotic covariance of the
maximum likelihood estimator $\ML(X_1,\dots,X_n).$ 
In the symmetric case $\Sigma=I_d$, it follows from \eqref{smed:asym} and \eqref{smed:asymcov:spez}
that \begin{equation}\label{eq:2-exp}
  \SMed(X_1,\dots,X_n) \; \dconv \; N_2(0,8I_2),
\end{equation}
which means that the asymptotic covariances of $\SMed(X_1,\dots,X_n)$ and 
$\ML(X_1,\dots,X_n)$ coincide in this case. This is not surprising, because      
$\SMed(X_1,\dots,X_n)$ is the maximum likelihood estimator of~$\mu$ 
in this symmetric case.

\smallbreak
(c) With $h(t)=(1+t)^{-3/2}$ we obtain a class of bivariate Cauchy distributions ${\rm MC}_2(\mu,\Sigma)$ with 
densities 
\begin{equation*}
    f(x;\mu,\Sigma)=\frac{1}{2\pi}\frac{1}{\sqrt{\det \Sigma}}
  \left (1 + (x-\mu)\t \Sigma^{-1}(x-\mu)\right )^{-\frac{3}{2}},~x\in\bR^2.
\end{equation*}
The associated marginal distributions are univariate Cauchy distributions. 
For this distribution family the sample mean is not even consistent, hence  we do not take it into
consideration as an estimator for $\mu$. Again we use \cite[formulas (3.2), (3.3)]{Mitchell} to see
that for given fixed positive definite $\Sigma$ the Fisher information
matrix is $\frac{3}{5}\Sigma^{-1}.$ Thus, by the general asymptotic theory 
of maximum likelihood estimators,
\begin{equation*}
   \ML(X_1,\dots,X_n) \; \dconv \; N_2\left (0,\tfrac{5}{3}\Sigma \right ).
\end{equation*}
Using
\begin{equation*}
   c_1(h)=2\int_0^\infty  (1+t^2)^{-\frac{3}{2}}\,dt
     =2,\hskip5mm  c_2(h)=2\pi \int_0^\infty t\hskip0.3mm (1+t^2)^{-\frac{3}{2}} \,dt=2\pi
\end{equation*}
we arrive at the value $\tfrac{\pi^2}{4}\Sigma$ for the asymptotic covariance matrix of the sample quarter median.
In the symmetric case $\Sigma=I_d$, using \eqref{smed:asym} and \eqref{smed:asymcov:spez}, it follows that
\begin{equation}\label{cauchy}
  \SMed(X_1,\dots,X_n) \; \dconv \; N_2(0,2I_2).
\end{equation}
  
\smallbreak
(d) With $h(t)=(1-t)^2$ for $0\le t \le 1$ and 0 elsewhere we have a special class of symmetric
bivariate Pearson type II distributions ${\rm SMPII}_2(\mu,\Sigma)$ with densities 
\begin{equation*}
   f(x;\mu,\Sigma)=\frac{3}{\pi}\frac{1}{\sqrt{\det \Sigma}}
        \max\bigl\{0, 1 - (x-\mu)\t \Sigma^{-1}(x-\mu)\bigr\}^2,~x\in\bR^2.
\end{equation*}
The asymptotic covariance matrix for the sample mean of independent two-dimensional
random vectors $X_1,\dots,X_n$ with the density $f(\cdot;\mu,\Sigma)$ is $\cov(X_1)=\frac{1}{8}\Sigma$; 
see \cite[p. 89]{Fang}.
The centered bivariate limit normal distribution of the maximum likelihood estimator has the covariance
matrix $\frac{1}{12}\Sigma$. Using
\begin{equation*}
    c_1(h)=2\int_0^1 (1-t^2)^2\,dt=\frac{16}{15},\hskip5mm  c_2(h)=2\pi \int_0^1 t\hskip0.3mm (1-t^2)^2\,dt=\frac{\pi}{3}
\end{equation*}
it follows that $\frac{25}{1024}\pi^2\Sigma$ is the asymptotic covariance matrix for the sample 
quarter median. In the symmetric case $\Sigma=I_d$ we obtain
\begin{equation*}
  \SMed(X_1,\dots,X_n) \; \dconv \; N_2\left (0,\frac{25}{128}I_2\right ).
\end{equation*}

\smallbreak
(e) Taking $h(t)=\frac{\exp(-t)}{\left (1+\exp(-t)\right )^2},\,t\ge 0,$ we obtain the symmetric bivariate logistic 
distributions ${\rm SML}_2(\mu,\Sigma)$, with densities given by
\begin{equation*}
  f(x;\mu,\Sigma)=\frac{1}{\pi\sqrt{\det \Sigma}}
  \frac{\exp\left (-(x-\mu)\t \Sigma^{-1}(x-\mu)\right )}
    {\left (1+\exp\left (-(x-\mu)\t \Sigma^{-1}(x-\mu)\right )\right )^2},~x\in\bR^2.
\end{equation*}
For $X\sim {\rm SML}_2(0,I_2)$ the distribution of $R^2=\|X\|^2$ is the 
univariate half-logistic distribution with density 
$2\frac{\exp(-t)}{\left (1+\exp(-t)\right )^2},\,t\ge 0$. Then, 
\begin{equation}\label{const1:mlogis}
 m_{{\rm SML}_2}:=\frac{1}{2}E(R^2)\,=\,\int_0^\infty t \frac{\exp(-t)}{\left (1+\exp(-t)\right )^2}\,dt\,\approx\,0.69314718
\end{equation}  
and $\cov(X)=m_{{\rm SML}_2}I_2.$ From this we deduce that the asymptotic covariance matrix
for the sample mean of independent two-dimensional random vectors $X_1,\dots,X_n$ with the 
above density
$f(\cdot;\mu,\Sigma)$ is $\cov(X_1)=m_{{\rm SML}_2}\Sigma$. 
With
\begin{equation}\label{const2:mlogis}
  k_{{\rm SML}_2}:=4\int_0^\infty t\left (\frac{1-\exp(-t)}{1+\exp(-t)}\right )^2\frac{\exp(-t)}
  {\left (1+\exp(-t)\right )^2}\,dt\,\approx\,1.59086291
\end{equation}
the Fisher information matrix is $k_{{\rm SML}_2}\Sigma^{-1}$, and the centered bivariate limit normal
distribution of the maximum likelihood estimator has the covariance matrix
$k_{{\rm SML}_2}^{-1}\Sigma$. Further, with 
\begin{equation*}
  c_1(h)=\int_0^\infty t^{-\tfrac{1}{2}}\frac{\exp(-t)}{\left (1+\exp(-t)\right )^2}\,dt
             \approx 0.67371824,\hskip3mm
  c_2(h)=\pi \int_0^\infty\frac{\exp(-t)}{\left (1+\exp(-t)\right )^2}\,dt=\frac{\pi}{2}
\end{equation*}
we obtain that $\frac{c_2(h)^2}{(2c_1(h))^2}\hskip0.5mm\Sigma\,
                    \approx\,1.35901156\hskip0.5mm\Sigma$
is the covariance matrix of the centered limiting normal distribution of the quarter median.   
In the symmetric case $\Sigma=I_d$ we get
\begin{equation*}
  \SMed(X_1,\dots,X_n) \; \dconv \; N_2\left (0,\frac{1}{2c_1(h)^2} I_2\right ),
\end{equation*}
where $\frac{1}{2c_1(h)^2}\approx 1.10157328$. Numerical values for the ratios $k_{{\rm SML}_2}^{-1}/m_{{\rm SML}_2}$
and $k_{{\rm SML}_2}^{-1}/(\frac{c_2(h)}{2c_1(h)})^2$ are shown in the last line of Table \ref{tab:eff}.
\end{example}

We summarize some special results presented in Example~\ref{ex:CLT} for the respective
strictly elliptical cases 
in Table~\ref{tab:eff}, where the relative asymptotic efficiencies 
${\rm eff}({\rm Mean}, {\rm ML})$, and ${\rm eff}(\QMed,{\rm ML})$ of the mean and the quarter median with respect to the maximum likelihood estimator are shown; these are the ratios
$\sigma_{\rm ML}^2/\sigma_{\rm Mean}^2$ and $\sigma_{\rm ML}^2/\sigma_{\QMed}^2$.

\begin{table}
  \caption{Relative asymptotic efficiencies}\label{tab:eff}
  \vskip4mm
  \hskip1.85cm
  \begin{tabular}{lcc}
    Distribution family          & ${\rm eff}({\rm Mean},{\rm ML})$  & ${\rm eff}(\QMed,{\rm ML})$  \\
    \noalign{\vspace{0.3mm}}
    \hline
    \noalign{\vspace{1mm}}
    bivariate normal             &         1                        &   $\frac{2}{\pi}$           \\
    \noalign{\vspace{1mm}}
    bivariate doubly exponential &        $\frac{2}{3}$             &   $\frac{8}{\pi^2}$                      \\
    \noalign{\vspace{1mm}}
    bivariate Cauchy             &         -                        &   $\frac{20}{3\pi^2}$                    \\
    \noalign{\vspace{1mm}}
    bivariate Pearson type II    &        $\frac{2}{3}$             &   $\frac{256}{75\pi^2}$     \\
    \noalign{\vspace{1mm}}
    symmetric bivariate logistic &        0.90686321                &   0.46253446    \\
  \end{tabular}
\end{table}

In contrast to the other multivariate versions of the median that appear
in the above example the quarter median goes beyond providing a location estimate
as it comes with an equipartition basis. In the elliptical case 
$P=\Ell_d(\mu,\Sigma;h)$ the basis elements are the unit vectors proportional to the 
eigenvectors of $\Sigma$ and thus contain information about the dispersion parameter. 
In Theorem~\ref{thm:CLTtheta} the asymptotic behavior of $(U_n)_{n\in\bN}$ is important.
In the strictly elliptical case we have consistency and we may assume that convergence takes
place in $\SO(d)$; see also the beginning of the proof of Theorem~\ref{thm:CLTtheta}. We aim 
at a second order statement about the distributional asymptotics of this sequence. 

In the following theorem we only consider $d=2$, but see the ensuing remarks. 
We can then use the parametrization of $\SO(2)$ given in Remark~\ref{rem:uniqueness}\,(b)
by an angle from an half-open interval $I$ of length $\pi/2$. We assume that $U=U_\alpha$ with
$I$ chosen such that $\alpha$ is in its interior. Then consistency implies that 
$U_n=U_{\alpha_n}$ with $I$-valued random variables $\alpha_n$ that converge to $\alpha$
with probability~1. It then makes sense to consider the distributional asymptotics of
$\sqrt{n}(\alpha_n-\alpha)$ as $n\to\infty$.

\begin{theorem}\label{thm:alpha}
Let $P=\Ell_2(\mu,\Sigma;h)$ and let $X_n, \theta_n, U_n$, $n\in\bN$, be as in
Theorem~\ref{thm:CLTell}\,\emph{(a)}. Let $\lambda_1$ and $\lambda_2$ be the eigenvalues of 
$\,\Sigma$, with $\lambda_1>\lambda_2$. Then, with $\alpha$ and $\alpha_n$
as defined above, 
\begin{equation*}
   \sqrt{n}(\alpha_n-\alpha)\; \dconv \; N\bigl(0,\sigma^2(\lambda_1,\lambda_2)\bigr),
\end{equation*}

\vspace{-2mm}\noindent
where 

\vspace{-5mm}
\begin{equation}\label{eq:asvarU}
    \sigma^2(\lambda_1,\lambda_2)\, =\, \frac{c_2(h)^2\lambda_1\lambda_2}
                                        {4c_1(h)^2(\sqrt{\lambda_1}-\sqrt{\lambda_2})^2}.
\end{equation}
Moreover, $\alpha_n$ and $\theta_n$ are asymptotically independent.
\end{theorem} 

In particular, the asymptotic variance in~\eqref{eq:asvarU} is large if the 
eigenvalues are close to each other. This is to be expected as  $\Ell_2(\mu,\Sigma;h)$
is then close to a symmetric distribution, where $U$ is not unique. 

In dimensions higher than two a possible approach could be based on the
representation of the  Lie group $\SO(d)$ by its Lie algebra $\so(d)$, which 
consists of the skew-symmetric matrices $A\in \bR^{d\times d}$, where
$U\in\SO(d)$ is written as the matrix exponential $U=\exp(A)$ of some
$A\in\so(d)$; see~\cite[Theorem A9.7]{Muirhead}~for a proof of the latter
assertion, and~\cite{Hall} for an elementary introduction to Lie groups, Lie algebras,
and representations. This leads to  an asymptotic
distribution $Q$ for $(U\t U_n)^{\sqrt{n}}$ as $n\to\infty$, where $Q$ is a probability measure 
on (the Borel subsets of) $\SO(d)$, and it avoids any ambiguities caused by parametrization. 
In fact, $Q$ is the distribution of the random orthogonal
matrix $\exp(A)$, where $A$ is a random skew symmetric matrix with jointly
normal subdiagonal entries. If $d=2$, as in the above theorem, then
\begin{equation*}
   U_\alpha = \exp\bigl(A(\alpha)\bigr)\quad\text{with } 
                    A:= \begin{pmatrix} 0 & -\alpha\\ \alpha & 0\end{pmatrix},
\end{equation*}
and $A$ is then specified by its one subdiagonal entry, which has a central normal 
distribution with variance given by~\eqref{eq:asvarU}.

\subsection{Measurability, selection and algorithms}\label{subsec:select}
Let $S$ be a topological space with Borel $\sigma$-field $\cB(S)$. We write $C_b(S)$ 
for the set of bounded continuous functions $f:S\to\bR$. Let $Z,Z_1,Z_2,\ldots$ be
$S$-valued random variables, i.e.~$\cB(S)$-measurable functions on some probability
space (which may depend on the respective variable). In its classical form, convergence
in distribution $Z_n\dconv Z$ of $Z_n$ to $Z$ as $n\to\infty$ means that
\begin{equation}\label{eq:dconvclassical}
 \lim_{n\to\infty}Ef(Z_n) = Ef(Z) \quad\text{for all }f\in C_b(S).
\end{equation}
This is the notion that we used in Section~\ref{subsec:asymp}.
An extension of this concept, due to Hoffmann-J{\o}rgensen, can be applied if 
the $Z_n$'s are not measurable with respect to the full Borel $\sigma$-field on $S$; 
roughly, the expectations $Ef(Z_n)$ in~\eqref{eq:dconvclassical} are then 
replaced by outer expectations $E^\star f(Z_n)$.
Similarly, almost sure convergence now refers to outer probabilities.
The research monographs of~\cite{Dudley} and~\cite{vdV+W} give an in-depth treatment
of this circle of ideas, together with a variety of applications. 
This extension appears in our proof of Theorem~\ref{thm:CLTtheta} in connection
with empirical processes; it can also be used if the estimators $(\theta_n,U_n)$ 
are not measurable. We refer the reader to the paper of~\cite{KuelbsZinn}, where this 
is worked out in detail for the related situation of quantile processes. 

While this avoids the problem of choosing an estimator from the respective solution
set in a measurable way, it is of independent interest whether such a selection 
is possible.  For the classical one-dimensional median and
the marginal median this can obviously be done by choosing the midpoint of the
respective (component) interval, but no such simple rule seems to exist for the
quarter median. 

We think of estimators as functions of the random variables
$X_i$, $i\in\bN$, and it is then enough to establish measurability for these functions. 
To be precise, for a given 
$n\in\bN$ and an $n$-tuple $(x_1,\dots,x_n)\in (\bR^d)^n$ we denote by 
$\big (\theta_n(x_1,\dots,x_n),U_n(x_1,\dots,x_n)\big )\in \bR^d\times \SO(d)$ 
a solution of the quarter median problem for the empirical distribution 
$P_{n;x_1,\dots,x_n}=n^{-1}\sum_{k=1}^n\delta_{x_k}$. 

\begin{proposition}\label{prop:mess}
There exists a function 
$\tau:(\bR^d)^n\rightarrow \bR^d\times \SO(d)$ that is measurable
with respect to the Borel $\sigma$-algebra on $(\bR^d)^n$ and the Borel $\sigma$-algebra
on $\bR^d\times \SO(d)$ such that, for each $(x_1,\dots,x_n)\in (\bR^d)^n$,
\begin{equation*}
\big (\theta_n(x_1,\dots,x_n),U_n(x_1,\dots,x_n)\big )\, :=\, \tau(x_1,\dots,x_n)
\end{equation*}
is a solution of the quarter median problem for $P_{n;x_1,\dots,x_n}$. 

Moreover,
$\tau$ may be chosen to be permutation invariant in the sense that 
\begin{equation}\label{eq:perminv}
  \tau(x_{\pi(1)},\ldots,x_{\pi(n)}) = \tau(x_1,\ldots,x_n)
\end{equation}
for all permutations $\pi$ of $\{1,\ldots,n\}$. 
\end{proposition}

Permutation invariance means that the estimates depend on the data only through the
respective empirical distribution.

\begin{figure}
\vspace{-.2cm}
\setlength{\abovecaptionskip}{-0.2cm}
\hskip0.65cm\includegraphics[scale=.5]{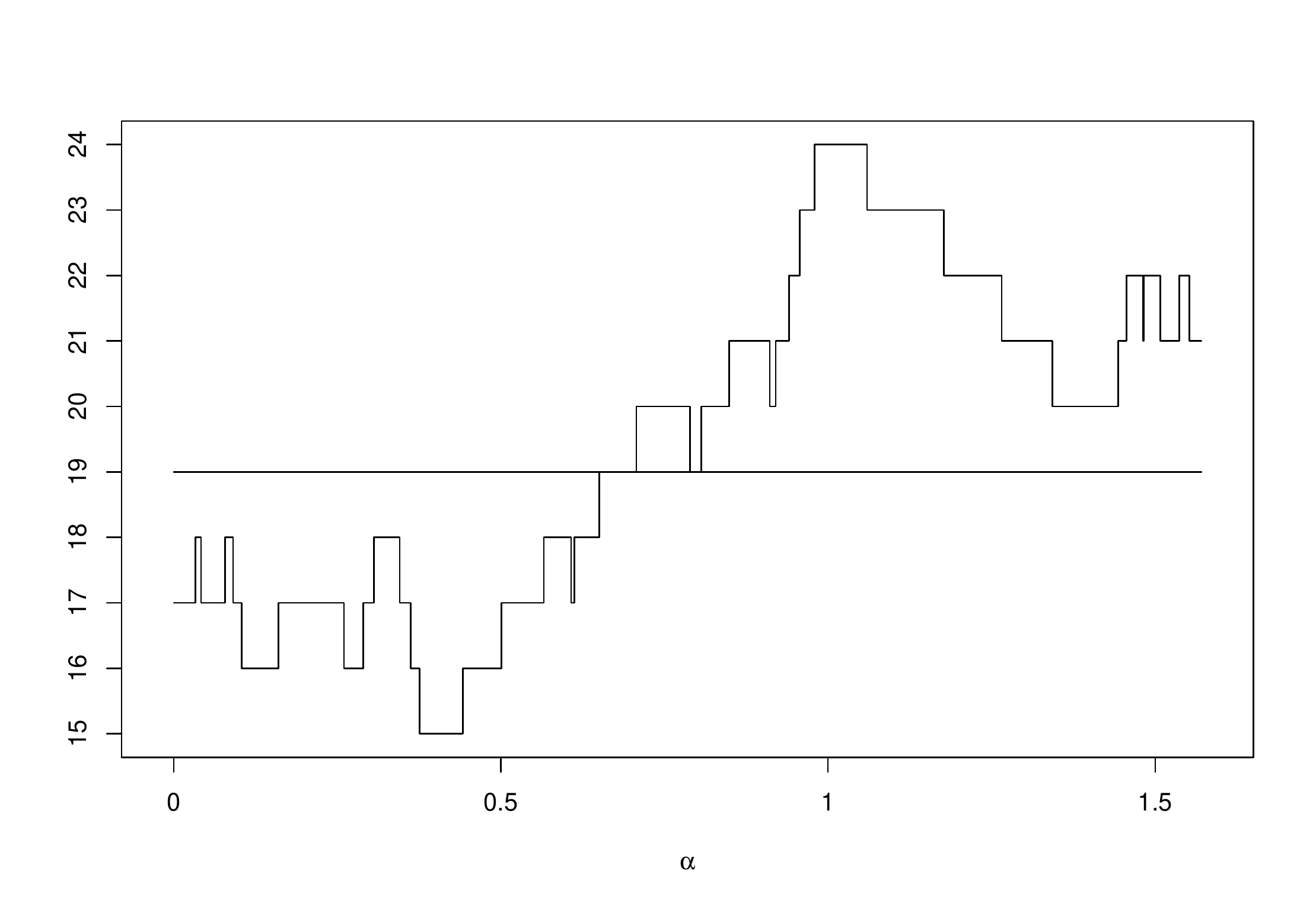}
\caption{Number of cities in the upper right quadrant, $\alpha\in[0,\pi/2)$}
\label{fig:staedtecount}
\end{figure}

Selecting a solution in a measurable way from the respective 
set of all solutions may seem as a corollary to establishing 
an algorithm that returns an estimate for every data
input $x_1,\ldots,x_n$. Indeed, without such an algorithm  
a statistical procedure would seem to be of limited use.
In dimension two,  and using the parametrization in Remark~\ref{rem:uniqueness}\,(b), 
we obtain a real function on the interval $[0,\pi/2)$ by counting the number of data 
points in the upper right quadrant specified by $U_\alpha$ and $\mMed_m(P^{U_\alpha})$; 
see Figure~\ref{fig:staedtecount} for a plot of this function for the city data
in Figure~\ref{fig:staedte}, where the counts refer to the open quarter space. 
By construction, all half-spaces contain at least half of
the data points, so each angle that leads to a count of $n/4$ (if $n$ is divisible by 4) 
would give a quarter median.
For absolutely continuous distributions this approach, together with the
intermediate value theorem, leads to a proof of Theorem~\ref{thm:main2} in dimension~2,
and this may be used numerically via bisection. Such a bisection argument
also works for functions that are piecewise constant and have jumps of size $\pm 1$ only,
provided a solution exists.

In what follows, still considering the case $d=2$, we present a different algorithm. 
Let $x_1,\dots,x_n$ be $n\ge 2$ pairwise distinct data points in $\bR^2$. 
For $1\le i<j\le n$ let $b_{ij}=\frac{x_j-x_i}{\|x_j-x_i\|}\in \mathbb{S}_1$
be the direction of the line through $x_i$ and $x_j,$ and let 
$b_{ij}'\in \mathbb{S}_1$ be orthogonal to $b_{ij}.$ 
With the empirical distribution $P_n=n^{-1}\sum_{k=1}^n \delta_{x_k}$ we associate 
the finite set $\mathscr{L}(P_n)$ of pairs $(\theta_{ij},U_{ij})$, with $U_{ij}\in \OO(2)$ 
as the matrix with row vectors $b_{ij}\t$ and ${b_{ij}'}\t$, and 
$\theta_{ij}=U_{ij}\t\eta$ with some $\eta\in\cM\bigl(P_n^{U_{ij}}\bigr)$.

\begin{theorem}\label{thm:algo} 
There is an element $(\theta_{ij},U_{ij}) \in \mathscr{L}(P_n)$
that solves the quarter median problem for $P_n$.
\end{theorem}

Hence, checking successively the conditions (5) and (6) for all
$\left (\theta_{ij},U_{ij}\right )\in \mathscr{L}(P_n)$, $1\le i<j\le n$,
provides a solution of the quarter median problem for $P_n$. This procedure,
however, may lead to an estimator that is not permutation invariant in the sense 
of~\eqref{eq:perminv}, meaning  that it would 
not be a function of the empirical distribution $P_n$. For our introductory data 
example, see Figure~\ref{fig:staedte}, this happens for both of the above algorithms:
Figure~\ref{fig:staedtecount} shows that different solutions appear depending on whether 
the angles are scanned clockwise or counterclockwise. Similarly, the algorithm based on 
Theorem~\ref{thm:algo} could lead to the $(i,j)$-pair corresponding to the cities 
Berlin and Frankfurt am Main, or München and Dortmund, for example. 

\section{Numerical comparison}
We present a small simulation study on the performance of the quarter median (\QMed) as 
a location estimator, comparing it with that of three other procedures---the spatial 
median (\SMed); Oja's simplicial median (\OMed), see \cite{OjaMed}; and 
Tukey's half-space median (\TMed), see \cite{Tukey}.
For a detailed description of these classical estimators we refer to the survey 
papers of~\cite{Sm1990} and \cite{Oja}. A Monte Carlo study on the performance of 
the three (and some other bivariate location) estimators
is given by \cite{MassPlant}.  

Samples are drawn from the strictly elliptical
distributions $N_2(\mu,\Sigma)$ (bivariate normal), ${\rm MDE}_2(\mu,\Sigma)$ (bivariate doubly exponential),
${\rm SMPII}_2(\mu,\Sigma)$ (symmetric bivariate Pearson type II), and ${\rm SML}_2(\mu,\Sigma)$
(symmetric bivariate logistic) introduced in Example~\ref{ex:CLT}.
The bivariate Cauchy distributions ${\rm MC}_2(\mu,\Sigma)$ also discussed
there are not taken into consideration for
the simulation study, as the Oja median does not exist for this distribution family. In fact,  
for a given distribution $P$ on $\bR^2$, independent $X_j\sim P$, $j=1,2$, and $\theta\in\bR^2$,
let $\Delta(X_1,X_2,\theta)$ be the area of the triangle in $\bR^2$ with vertices $X_1,X_2,\theta.$ 
If the expectation $\gamma_P(\theta):=E\left (\Delta(X_1,X_2,\theta)\right )$ is finite for all $\theta\in\bR^2$, the
Oja median is a point $\theta_{\rm Oja}(P)\in\bR^2$ that minimizes the function
$\theta\mapsto \gamma_P(\theta)$. Obviously, $\gamma_P\equiv\infty$ for $P={\rm MC}_2(0,\Sigma)$, while
for the other elliptical distributions $P$ under consideration
$\theta_{\rm Oja}(P)=\mu.$ 

Without loss of generality we take $\mu=0$. The dispersion matrix chosen is
$\Sigma=\begin{pmatrix} 1 & 0 \\ 0 & \lambda \end{pmatrix}$, with
$\lambda\in \{0.01,0.1,0.5,0.9\}$. Repeatedly, with $r=10000$ replications, samples of size $n=100$ are drawn 
from the underlying distribution. Let $\SMed_{n,i},\OMed_{n,i},\TMed_{n,i},\QMed_{n,i}$
be the observed values of the corresponding estimators obtained in the $i$th repetition.    
Table \ref{num:sim} shows the numerical values of the components $\hat m_1,\hat m_2$ of the empirical means and
the numerical values of the eigenvalues $\hat l_1,\hat l_2$ of the empirical covariance matrices
of $\sqrt{n}\SMed_{n,i},\sqrt{n}\OMed_{n,i},\sqrt{n}\TMed_{n,i},\sqrt{n}\QMed_{n,i}$, $i=1,\ldots,r.$  
The simulations are conducted by using the statistical software environment {\bf R};
the calculation of the estimators \SMed, \OMed,~and \TMed~is easily done by applying the corresponding functions 
{\tt med(\ldots,method="Spatial"), med(\ldots,method="Oja",\ldots)}~and {\tt TukeyMedian(\ldots)} 
provided with the additional {\bf R} software packages {\tt depth} and {\tt TukeyRegion}.
Note that, in the later package, the Tukey median is defined to be the barycenter of the 
Tukey region.

\hskip5mm
\begin{longtable}{rrrrrr}
   \caption{Components $\hat m_1,\hat m_2$ of the empirical means, and eigenvalues $\hat l_1,\hat l_2$ of the empirical
     covariance matrices of the $\sqrt{n}$-scaled estimators~\SMed, \OMed, \TMed, \QMed~obtained by simulations with 10000
     replications for strictly elliptical distributions with dispersion matrix
     $\Sigma=\protect\begin{pmatrix} 1 & 0 \\ 0 & \lambda \protect \end{pmatrix}$,
     $\lambda\in \{0.01,0.1,0.5,0.9\}$; sample size $n=100.$}
     \label{num:sim}\\

    \noalign{\vspace{8mm}}    

    \multicolumn{6}{c}{$\lambda=0.01$} \\

    \noalign{\vspace{3mm}}
    & &\multicolumn{1}{c}{$\hat m_1$} & \multicolumn{1}{c}{$\hat m_2$} & \multicolumn{1}{c}{$\hat l_1$} & \multicolumn{1}{c}{$\hat l_2$} \\

    \noalign{\vspace{1mm}}
    \hline

    \noalign{\vspace{2mm}}
    \multirow{4}{*}{$N_2(0,\Sigma)$} & \SMed & 0.009648 & -0.002215 & 1.453955 & 0.019903 \\
    & \OMed & 0.001726 & -0.001769 & 1.282321 & 0.012935 \\
    & \TMed &-0.000700 & -0.001661 & 1.311815 & 0.013135 \\
    & \QMed & 0.010454 & -0.001571 & 1.558638 & 0.015973 \\

    \noalign{\vspace{2mm}}
    \multirow{4}{*}{${\rm MDE}_2(0,\Sigma)$} & \SMed & 0.062843 & 0.000282 & 9.651512 & 0.132048 \\
    & \OMed &  0.047608 & -0.001115 & 8.376979 & 0.084180 \\
    & \TMed &  0.041778 & -0.001260 & 8.548116 & 0.085507 \\
    & \QMed &  0.066367 &  0.000427 &10.363242 & 0.102538 \\

    \noalign{\vspace{2mm}}
    \multirow{4}{*}{${\rm SMPII}_2(0,\Sigma)$} & \SMed & -0.002203 & 0.000669 & 0.227956 & 0.003093 \\
    & \OMed & -0.000128 & 0.000219 & 0.201077 & 0.001999 \\
    & \TMed &  0.000433 & 0.000230 & 0.205588 & 0.002041 \\
    & \QMed & -0.003445 & 0.000139 & 0.243375 & 0.002450 \\
    
    \noalign{\vspace{2mm}}
    \multirow{4}{*}{${\rm SML}_2(0,\Sigma)$} & \SMed & -0.004630 & -0.000134 & 1.245813 & 0.017274 \\
    & \OMed & -0.014100 & -0.000211 & 1.102237 & 0.010843 \\
    & \TMed & -0.014277 & -0.000292 & 1.127860 & 0.011012 \\
    & \QMed & -0.001921 & -0.000723 & 1.326942 & 0.013266 \\
    
    \noalign{\vspace{8mm}}
    \multicolumn{6}{c}{$\lambda=0.1$} \\

    \noalign{\vspace{3mm}}
    & &\multicolumn{1}{c}{$\hat m_1$} & \multicolumn{1}{c}{$\hat m_2$} & \multicolumn{1}{c}{$\hat l_1$} & \multicolumn{1}{c}{$\hat l_2$} \\

    \noalign{\vspace{1mm}}
    \hline

    \noalign{\vspace{2mm}}
    \multirow{4}{*}{$N_2(0,\Sigma)$} & \SMed & 0.003782 & -0.002045 & 1.347796 & 0.139426 \\
    & \OMed &  0.002807 & -0.000513 & 1.290243 & 0.127393 \\
    & \TMed &  0.004400 & -0.001163 & 1.316496 & 0.129837 \\
    & \QMed &  0.001002 & -0.002581 & 1.550401 & 0.157899 \\

    \noalign{\vspace{2mm}}
    \multirow{4}{*}{${\rm MDE}_2(0,\Sigma)$} & \SMed & -0.024886 & -0.014379 & 8.864224 & 0.937143 \\
    & \OMed & -0.025370 & -0.011321 & 8.434326 & 0.854891 \\
    & \TMed & -0.020137 & -0.009642 & 8.564605 & 0.873420 \\
    & \QMed & -0.031408 & -0.003320 & 10.250833 & 1.052323 \\

    \noalign{\vspace{2mm}}
    \multirow{4}{*}{${\rm SMPII}_2(0,\Sigma)$} & \SMed & -0.008639 & -0.002740 & 0.204793 & 0.021491 \\
    & \OMed & -0.006975 & -0.002715 & 0.197094 & 0.019512 \\
    & \TMed & -0.006893 & -0.002416 & 0.201236 & 0.019897 \\
    & \QMed & -0.011107 & -0.001635 & 0.232204 & 0.023828 \\

    \noalign{\vspace{2mm}}
    \multirow{4}{*}{${\rm SML}_2(0,\Sigma)$} & \SMed & -0.008358 & -0.006214 & 1.175185 & 0.120428 \\
    & \OMed & -0.011263 & -0.004974 & 1.119404 & 0.110270 \\
    & \TMed & -0.013367 & -0.005026 & 1.146522 & 0.112307 \\
    & \QMed & -0.007513 & -0.001273 & 1.351009 & 0.136288 \\
    
    \noalign{\vspace{8mm}}
    \multicolumn{6}{c}{$\lambda=0.5$} \\

    \noalign{\vspace{3mm}}
    & &\multicolumn{1}{c}{$\hat m_1$} & \multicolumn{1}{c}{$\hat m_2$} & \multicolumn{1}{c}{$\hat l_1$} & \multicolumn{1}{c}{$\hat l_2$} \\

    \noalign{\vspace{1mm}}
    \hline

    \noalign{\vspace{2mm}}
    \multirow{4}{*}{$N_2(0,\Sigma)$} & \SMed & 0.011107 & -0.007306 & 1.274190 & 0.639970 \\
    & \OMed & 0.011333 & -0.006804 & 1.274619 & 0.636957 \\
    & \TMed & 0.011808 & -0.005540 & 1.309882 & 0.650354 \\
    & \QMed & 0.009552 & -0.007997 & 1.551297 & 0.767945 \\

    \noalign{\vspace{2mm}}
    \multirow{4}{*}{${\rm MDE}_2(0,\Sigma)$} & \SMed & -0.022503 & -0.017123 & 8.333408 & 4.103622 \\
    & \OMed & -0.019896 & -0.016413 & 8.340335 & 4.101557 \\
    & \TMed & -0.018584 & -0.020229 & 8.551098 & 4.201095 \\
    & \QMed & -0.022966 & -0.005561 & 10.177384 & 5.088397 \\

    \noalign{\vspace{2mm}}
    \multirow{4}{*}{${\rm SMPII}_2(0,\Sigma)$} & \SMed & -0.005249 & -0.001641 & 0.197165 & 0.097735 \\
    & \OMed & -0.005218 & -0.001727 & 0.196625 & 0.097539 \\
    & \TMed & -0.005330 & -0.002090 & 0.200822 & 0.099652 \\ 
    & \QMed & -0.006836 & -0.000873 & 0.237224 & 0.118714 \\

    \noalign{\vspace{2mm}}
    \multirow{4}{*}{${\rm SML}_2(0,\Sigma)$} & \SMed & 0.008670 & 0.013723 & 1.099552 & 0.563826 \\
    & \OMed & 0.009137 & 0.013315 & 1.101891 & 0.561546 \\
    & \TMed & 0.006853 & 0.015428 & 1.125076 & 0.573525 \\
    & \QMed & 0.012909 & 0.018256 & 1.337041 & 0.684311 \\
                                                   
    \noalign{\vspace{12mm}}
    \multicolumn{6}{c}{$\lambda=0.9$} \\

    \noalign{\vspace{3mm}}
    & &\multicolumn{1}{c}{$\hat m_1$} & \multicolumn{1}{c}{$\hat m_2$} & \multicolumn{1}{c}{$\hat l_1$} & \multicolumn{1}{c}{$\hat l_2$} \\

    \noalign{\vspace{1mm}}
    \hline

    \noalign{\vspace{2mm}}
    \multirow{4}{*}{$N_2(0,\Sigma)$} & \SMed & 0.010573 & 0.002084 & 1.284202 & 1.157361 \\
    & \OMed & 0.010942 & 0.001677 & 1.291926 & 1.166358 \\
    & \TMed & 0.009962 & 0.002722 & 1.314541 & 1.191291 \\
    & \QMed & 0.007804 & 0.005679 & 1.553846 & 1.414504  \\

    \noalign{\vspace{2mm}}
    \multirow{4}{*}{${\rm MDE}_2(0,\Sigma)$} & \SMed & 0.036521 & 0.070980 & 8.360825 & 7.417609 \\
    & \OMed & 0.037581 & 0.072080 & 8.417271 & 7.462682 \\
    & \TMed & 0.037717 & 0.070611 & 8.586432 & 7.609654 \\
    & \QMed & 0.007530 & 0.087122 &10.274479 & 9.154238 \\

    \noalign{\vspace{2mm}}
    \multirow{4}{*}{${\rm SMPII}_2(0,\Sigma)$} & \SMed & -0.006974 & -0.000672 & 0.193231 & 0.178529 \\
    & \OMed & -0.006773 & -0.001132 & 0.194627 & 0.178877 \\
    & \TMed & -0.006907 & -0.001222 & 0.198142 & 0.183445 \\ 
    & \QMed & -0.005604 & -0.002042 & 0.234378 & 0.220361 \\

    \noalign{\vspace{2mm}}
    \multirow{4}{*}{${\rm SML}_2(0,\Sigma)$} & \SMed & -0.003381 & 0.003081 & 1.078231 & 0.982265 \\
    & \OMed & -0.005205 & 0.002348 & 1.082106 & 0.989597 \\
    & \TMed & -0.002501 & 0.003175 & 1.107970 & 1.007464 \\
    & \QMed & -0.010836 &-0.000977 & 1.321225 & 1.190672  
  \end{longtable}
For samples taken from bivariate normal distributions with dispersion matrix
$\Sigma=\begin{pmatrix} 1 & 0 \\ 0 & \lambda \end{pmatrix}$ the limit distributions of
\SMed, \OMed, and \QMed~as $n\to\infty$ are the centered bivariate normal distributions with
the diagonal covariance matrices 
$\Sigma_{\SMed}=\begin{pmatrix} l_{1,\SMed} & 0 \\ 0 & l_{2,\SMed} \end{pmatrix}$,
see \cite{BrownSMed}, where also expressions for the diagonal elements 
$l_{1,\SMed},l_{2,\SMed}$ and their numerical
values for some special values of $\lambda$ are given, 
$\Sigma_{\OMed}=\begin{pmatrix} 4/\pi & 0 \\ 0 & 4\lambda/\pi \end{pmatrix}$, see \cite{OjaNii}, and
$\Sigma_{\QMed}=\begin{pmatrix} \pi/2 & 0 \\ 0 & \pi\lambda/2 \end{pmatrix}$.
The limit distribution of $\TMed$ is the distribution of
$\arg \sup_{t\in\bR^2} \inf_{u\in \Seins} \left ( Z(u) - \frac{1}{2\pi \lambda^{1/2}} u\t t\right )$,
where $\left (Z(u),u\in \Seins\right )$ is a centered Gaussian process with the covariance function
$\cov\left (Z(u),Z(v)\right )=P(V_{++}(0;u,v)-1/4,\,u,v\in\Seins;$ see \cite{Nolan}. To the best of
our knowledge, nothing is known about the type of this distribution. In this respect, see for example also
the interesting results on the asymptotic behaviour of the empirical Tukey depth process given by
\cite{Mass2004}.

Table \ref{emp:lim:comp} shows the
values of $\hat l_1,\hat l_2$ and, for comparison, that of the diagonal elements $l_1,l_2$ of
$\Sigma_{\SMed},\Sigma_{\OMed}$, and $\Sigma_{\QMed}$. 

\begin{table}[ht]
  \begin{center}
    \caption{Bivariate normal case: Empirical eigenvalues $\hat l_1,\hat l_2$ and eigenvalues $l_1,l_2$ of
 $\Sigma_{\SMed},\Sigma_{\OMed}$, and $\Sigma_{\QMed}$}\label{emp:lim:comp} 
  \begin{tabular}{rrrrrr}
    &       & \multicolumn{1}{c}{$\hat l_1$} & \multicolumn{1}{c}{$\hat l_2$} & \multicolumn{1}{c}{$l_1$} &\multicolumn{1}{c}{$l_2$} \\
    \noalign{\vspace{1mm}}
    \hline
   \noalign{\vspace{1mm}}
   \multirow{3}{*}{$\lambda=0.01$} & \SMed & 1.453955 & 0.019903 & 1.555954 & 0.000624 \\ 
    &\OMed & 1.282321 & 0.012935 & 1.273240 & 0.012732 \\
    &\QMed & 1.558638 & 0.015973 & 1.570796 & 0.015708 \\

    \noalign{\vspace{2mm}}
    \multirow{3}{*}{$\lambda=0.1$} & \SMed & 1.347796 & 0.139426 & 1.460156 & 0.019491 \\ 
    &\OMed &  1.290243 & 0.127393 & 1.273240 & 0.127324 \\
    &\QMed &  1.550401 & 0.157899 & 1.570796 & 0.157080 \\

    \noalign{\vspace{2mm}}
    \multirow{3}{*}{$\lambda=0.5$} & \SMed & 1.274190 & 0.639970 & 1.304688 &  0.329485 \\ 
    &\OMed & 1.274619 & 0.636957 & 1.273240 & 0.636620 \\
    &\QMed & 1.551297 & 0.767945 & 1.570796 & 0.785398 \\

   \noalign{\vspace{2mm}}
   \multirow{3}{*}{$\lambda=0.9$} & \SMed & 1.284202 & 1.157361 & 1.274100 & 1.032058 \\ 
    &\OMed & 1.291926 & 1.166358 & 1.273240 & 1.145916 \\
    &\QMed & 1.553846 & 1.414504 & 1.570796 & 1.413717       
  \end{tabular}
  \end{center}
\end{table}

\section{Proofs}\label{sec:proofs} 

\subsection{Proof of Theorem~\ref{thm:main2}}\label{subsec:thm1}
Using Makeev's result, we only need to remove the smoothness assumption.

For an arbitrary probability measure $P$ on $(\bR^d,\cB^d)$ let 
$P_\epsilon=P * N_d(0,\epsilon I_d)$ be the convolution of $P$ with a centered 
$d$-dimensional normal distribution with independent components that all have variance
$\epsilon>0$. Then $P_\epsilon$ is absolutely continuous with an everywhere positive density, 
so that a solution $\bigl(\theta_\epsilon, U_\epsilon)$  
of the equipartition problem exists. Let $\vvvert \cdot \vvvert$ be a matrix norm
on $\bR^{d\times d}$ that is compatible with the Euclidean norm $\| \cdot \|$ on $\bR^d$. 
From $\theta_\epsilon=U_\epsilon\t U_\epsilon \theta_\epsilon$ it then follows that
\begin{equation}\label{ineq:qmed}
  \|\theta_\epsilon\|\le \vvvert U_\epsilon\t \vvvert \|U_\epsilon \theta_\epsilon \| \le s \|U_\epsilon \theta_\epsilon \|,
\end{equation}
where $s=\sup_{U\in \SO(d)}\vvvert U\t \vvvert<\infty.$ Let
$t=\limsup_{\epsilon\to 0} \| U_\epsilon\theta_\epsilon \|$. Hence there
exists a sequence $(\epsilon_n)_{n\in\bN}$ with $\epsilon_n\to 0$ such that
$\| U_{\epsilon_n}\theta_{\epsilon_n}\| \rightarrow t $ and $U_{\epsilon_n}\to U \in\OO(d)$ 
as $n\to\infty$.
As a consequence, $P_{\epsilon_n}^{U_{\epsilon_n}}\rightarrow P^U$ weakly. Due to the fact that $U_{\epsilon_n}\theta_{\epsilon_n}$
is a marginal median of $P_{\epsilon_n}^{U_{\epsilon_n}}$ we deduce from this that the sequence
$\left (\| U_{\epsilon_n}\theta_{\epsilon_n}\|\right )_{n=1}^\infty$ is bounded; hence, 
by~\eqref{ineq:qmed}, the sequence $(\theta_{\epsilon_n})_{n=1}^\infty$ is also bounded. 
Therefore, $\theta_{\epsilon_{n'}}\to\theta\in \bR^d$ and
$U_{\epsilon_{n'}}\to U \in\OO(d)$ along some subsequence $(\epsilon_{n'})$ of $(\epsilon_{n}).$
To see that
$(\theta,U)$ solves the quarter median problem for $P$ we use that $P_{\epsilon_{n'}}^{U_{\epsilon_{n'}}}\rightarrow P^U$
weakly and argue as follows: With
$b_{\epsilon_{n'},1}\t,\dots,b_{\epsilon_{n'},d}\t$ and $b_1\t,\dots,b_d\t$
as the row vectors of $U_{\epsilon_{n'}}$ and $U,$ respectively, we have that, for each $1\le i<j\le d$ and each $\delta>0$,
\begin{align*}\frac{1}{4}&\ \le \ \limsup_{n'\to\infty}P_{\epsilon_{n'}}\big (V_{+-}(\theta_{\epsilon_{n'}};b_{\epsilon_{n'},i},b_{\epsilon_{n'},j})\big )\\
                         &\ \le \ \limsup_{n'\to\infty}P_{\epsilon_{n'}}\big (b_{\epsilon_{n'},i}\t x \ge b_{i}\t \theta -\delta, \, b_{\epsilon_{n'},j}\t x \le b_{i}\t \theta + \delta \big )\\
                         &\ \le \ P\big (b_{i}\t x \ge  b_{i}\t \theta -\delta, \, b_{j}\t x \le  b_{i}\t \theta + \delta\big).
\end{align*}
Let $\delta\downarrow 0$ to obtain $P\big (V_{+-}(\theta;b_i,b_j)\big )\ge \frac{1}{4}.$ 
The other inequalities stated in \eqref{eq:qmeddef1} and \eqref{eq:qmeddef2} can be verified in the same
way.

\subsection{Proof of Proposition~\ref{prop:qumedequi}}

(a) Let $\theta$ be a quarter median for $P$ with $b_1,\ldots,b_d$ an associated 
orthonormal basis 
of $\bR^d$ such that the constraints~\eqref{eq:qmeddef1} and~\eqref{eq:qmeddef2} are
satisfied. Further, let $S(x)=Ax+c$ with $A\in\OO(d)$, $c\in\bR^d$, be a Euclidean motion. 
Then $b_1',\ldots,b_d'$ with $b_i':=Ab_i$, $1\le i\le d$, is again an orthonormal basis
for $\bR^d$, and it is easily checked that
\begin{align*}
  &S^{-1}\bigl(H_{\pm}(A\theta + c;b_i')\bigr) 
                       \, = \,  H_{\pm}(\theta;b_i),\quad 1\le i \le d,\\
  &S^{-1}\bigl(V_{\pm\pm}(A\theta + c;b_i',b_j')\bigr) 
                       \, = \,  V_{\pm\pm}(\theta;b_i,b_j),\quad 1\le i<j\le d.
\end{align*}
In view of the definition of push-forwards this implies that $S(\theta)$
is a quarter median for $P^S$, with $b_1',\ldots,b_d'$ an associated equipartition basis. 

\vspace{1mm}
(b) If $V=AU$ with $A\in G$ then
\begin{equation*}
   \Psi(V,P)\, =\, (AU)\t\mMed_m\bigl((P^U)^A\bigr)\, 
              =\, U\t A\t A \mMed_m(P^U)\,
              =\, \Psi(U,P),
\end{equation*}
where we have used~\eqref{eq:Ginv} and $G\subset \OO(d)$.

\vspace{1mm}
(c) Let $b_1\t,\ldots,b_d\t$ be the row vectors of $U$. With
\begin{equation*}
  \theta_U=(\theta_{U,1},\ldots,\theta_{U,d})\t\,:= \,
               U \theta= (b_1\t \theta,\ldots,b_d^t\theta)\t
\end{equation*}
it holds that 
$\theta=U\t\theta_U=\sum_{i=1}^d\theta_{U,i}b_i$ and that, by definition, $\theta_U$ is a  
marginal median of $P^U.$ We can now choose an element
$\Tilde\theta_U=(\Tilde\theta_{U,1},\ldots,\Tilde \theta_{U,d})\t \in \cM\left (P^U\right )$ such that
\begin{align*}
  P(b_i\t x \ge \theta_{U,i})=P(b_i\t x \ge \Tilde\theta_{U,i}),\hskip5mm  P(b_i\t x \le \theta_{U,i})=P(b_i\t x \le \Tilde\theta_{U,i})
 \end{align*}
for $1\le i\le d$, and
\begin{align*}
  P(b_i\t x \ge \theta_{U,i},\, b_j\t x \ge \theta_{U,j})&=P(b_i\t x \ge \Tilde\theta_{U,i},\, b_j\t x \ge \Tilde \theta_{U,j})\\
  P(b_i\t x \ge \theta_{U,i},\, b_j\t x \le \theta_{U,j})&=P(b_i\t x \ge \Tilde\theta_{U,i},\, b_j\t x \le \Tilde \theta_{U,j})\\
  P(b_i\t x \le \theta_{U,i},\, b_j\t x \ge \theta_{U,j})&=P(b_i\t x \le \Tilde\theta_{U,i},\, b_j\t x \ge \Tilde \theta_{U,j})\\
  P(b_i\t x \le \theta_{U,i},\, b_j\t x \le \theta_{U,j})&=P(b_i\t x \le \Tilde\theta_{U,i},\, b_j\t x \le \Tilde \theta_{U,j})
\end{align*}
for $1\le i<j\le d$. 
Then putting $\Tilde\theta=U\t\Tilde\theta_U=\sum_{i=1}^d \Tilde\theta_{U,i} b_i$ we have
\begin{align*}
  P(b_i\t x \ge b_i\t \theta)=P(b_i\t x \ge b_i\t \Tilde\theta),\hskip5mm  P(b_i\t x \le b_i\t \theta)=P(b_i\t x \le b_i\t \Tilde\theta)
 \end{align*}
for $1\le i\le d$, and
\begin{align*}
  P\bigl (b_i\t x \ge b_i\t \theta,\, b_j\t x \ge b_j\t \theta\bigr )&=P\bigl (b_i\t x \ge b_i\t\Tilde\theta,\, b_j\t x \ge b_j\t\Tilde\theta\bigr )\\
  P\bigl (b_i\t x \ge b_i\t \theta,\, b_j\t x \le b_j\t \theta\bigr )&=P\bigl (b_i\t x \ge b_i\t\Tilde\theta,\, b_j\t x \le b_j\t\Tilde\theta\bigr )\\
  P\bigl (b_i\t x \le b_i\t \theta,\, b_j\t x \ge b_j\t \theta\bigr )&=P\bigl (b_i\t x \le b_i\t\Tilde\theta,\, b_j\t x \ge b_j\t\Tilde\theta\bigr )\\
  P\bigl (b_i\t x \le b_i\t \theta,\, b_j\t x \le b_j\t \theta\bigr )&=P\bigl (b_i\t x \le b_i\t\Tilde\theta,\, b_j\t x \le b_j\t\Tilde\theta\bigr )
\end{align*}
for $1\le i<j\le d$. Thus, $(\Tilde\theta,U)$ solves the quarter median problem for $P$. 

\vspace{1mm}
(d) The first statement may be regarded as an equivariance property of solution
pairs for the quarter median problem; its proof proceeds as in (a). For the second
statement we note that $A\in\OO(d)$ and then calculate, using the shift equivariance 
of the marginal median,
\begin{align*}
   \Psi(U\! A\t,P^S)\ &=\ (U\! A\t)\t \mMed_m\bigl((P^S)^{U\! A\t}\bigr) \\
                      &=\ A U\t\bigl(\mMed_m(P^{UA\t A}) + U\!A\t b)\\
                      &=\ A\Psi(U,P) + b\quad =\ S\bigl(\Psi(U,P)\bigr).
\end{align*}

(e) Let $U\in\OO(d)$ and $(U_n)_{n=1}^\infty$ be a sequence of elements in $\OO(d)$
converging to $U.$ Due to $P^{U_n}\rightarrow P^U$ weakly and the fact that
$P^{b}$ has a unique median for all $b\in S_d,$ it follows that
$\mMed_m(P^{U_n})\rightarrow \mMed_m(P^U).$ Then, from 
\begin{align*}
      \vvvert U_n&\t\mMed_m(P^{U_n}) - U\t\mMed_m(P^U) \vvvert\\
    &\hskip3mm = ~\vvvert U_n\t\mMed_m(P^{U_n}) - U_n\t\mMed_m(P^U)  +  U_n\t\mMed_m(P^U) - U\t\mMed_m(P^U) \vvvert\hspace{5mm}\\
    &\hskip3mm \le ~\vvvert U_n\t\vvvert \|\mMed_m(P^{U_n}) - \mMed_m(P^U)\|  + \vvvert U_n\t - U\t \vvvert \|\mMed_m(P^U)\|
\end{align*} 
and the compactness of $\OO(d)$ we deduce that $ U_n\t\mMed_m(P^{U_n})\rightarrow U\t\mMed_m(P^U).$ 

\subsection{Proof of Theorem~\ref{thm:consistency}}\label{subsec:proofconsistency}
We will make use of the fact that if $m_n,n\in\bN,$ are
medians of univariate distributions that converge weakly as $n\to\infty$ to a 
univariate distribution 
with unique median $m$, then $\lim_{n\to\infty} m_n=m$. Further, we only prove the second
part of the theorem as the arguments for (a) are similar.
  
Let $(\Omega,\cA,\bP)$ be a probability space on which the random vectors $X_1,X_2,\ldots$ 
are defined, so that $\bP^{X_j}=P$ for all $j\in\bN$. 
Then there is a $\bP$-null set $N\in\cA$ such that for each $\omega \in N^c$  
\begin{equation*}
   P_{n;\omega}:=\frac{1}{n}\sum_{j=1}^n\delta_{X_j(\omega)}\rightarrow P\hskip3mm\text{weakly}.
\end{equation*}
Fix $\omega\in N^c$; until further notice we omit the argument $\omega$ below.
By definition, the quarter median $\theta_n$ of $P_n$ can be written as $\theta_n=U_n\t\eta_n$ with some $U_n\in \OO(d)$ and  some
$\eta_n\in[\mMed_-(P_n^{U_n}),\mMed_+(P_n^{U_n})]$
(both may depend on $\omega\in N^c$). 
Let $(U_{n'})$ be a subsequence of $(U_n)$. As $\OO(d)$ is compact a subsubsequence
$(U_{n''})$ converges to some $V\in\OO(d)$, and we have  weak convergence
$P_{n''}^{U_{n''}}\to P^V$ by the extended continuous mapping theorem; 
see~\cite[Theorem 5.5]{BillCoPM}. As the distribution $P^b$ has a unique median for all $b\in \Sd$, 
the marginal median $\eta_n''$ of $P_{n''}^{U_{n''}}$ converges to the marginal median $\eta$ of $P^V$.
Now fix $i,j$ with $1\le i < j\le d$. Let $a_n'',b_n''$, $a,b$ be the corresponding rows 
in $U_n''$ and $V$ respectively.
Using the same arguments as in the proof of Theorem~\ref{thm:main2}
we then obtain for each $\epsilon>0$, and with $m$ instead of $n''$, 
\begin{align*}
  \frac{1}{4}\ &\le\ \limsup_{m} P_{m}\bigl( a_{m}\t x \le a_m\t \theta_m, \
               b_m\t x \ge b_m\t \theta_m \bigr)\\
       & \le\ \limsup_{m} P_{m}\bigl( a_m\t x \le a\t \eta + \epsilon, \
               b_m\t x \ge b\t \eta -\epsilon\bigr)\\
       & \le\ P\bigl(a\t x \le a\t \eta + \epsilon, \
               b\t x \ge b\t \eta -\epsilon\bigr),
\end{align*}
where in the first line we have used that $(\theta_m,U_m)$ solves the quarter median 
problem for $P_m$.
Letting $\epsilon\downarrow 0$ we get 
$P\left( V_{-+}(\eta; a,b) \right)\, \ge\, 1/4$.
More generally, and using the same argument, we obtain 
$P\bigl(V_{\pm\pm}(\eta;a,b)\bigr)\ge 1/4$ as well as  
$P\bigl(H_{\pm}(\eta;a)\bigr)\ge 1/2$ and $P\big(H_{\pm}(\eta(\omega);b)\bigr))\ge 1/2$. 
Recall that some of the quantities depend on the $\omega$ chosen above so that, 
for example,
\begin{equation*}
 P\bigl(H_+(\eta;a)\bigr)\; =\; P\bigl(H_+(\eta(\omega);a(\omega))\bigr)
   \; = \; 
    P\bigl(\{x\in\bR^d:\, \langle a(\omega),x-\eta(\omega)\rangle\ge 0\}\bigr).
\end{equation*}

In summary, we have shown that every limit point $(\eta(\omega),V(\omega))$
of the sequence $(\theta_n(\omega),U_n(\omega))_{n\in\bN}$ is a solution of the quarter
median problem for $P$.
Hence, if the solution is unique and given by $(\theta,H)\in\bR^d\times \bH$, 
then, on a set of probability~1, $\eta=\theta$ and $[V]=H$.

\subsection{Proof of Theorem~\ref{thm:CLTtheta}}\label{subsec:proofCLTtheta}
Invariance with respect to shifts means that we may assume $\theta=0$. 
As already noticed in Section~\ref{subsec:proofconsistency}, the quarter median $\theta_n$ can be written as $\theta_n=U_n\t \eta_n$ with 
$\eta_n\in [\mMed_-(P_n^{U_n}),\mMed_+(P_n^{U_n})]$.
Theorem~\ref{thm:consistency} implies that $([U_n])_{n\in\bN}$ converges almost 
surely to $H$  
with respect to  the quotient topology on $\bH(d)$. The mapping $p:\OO(d)\to \bH(d)$,
$U\mapsto [U]$, associates to each $H\in\bH(d)$ a finite fiber $p^{-1}(\{H\})$, and $\OO(d)$
is a covering space for $\bH(d)$. General results from algebraic 
topology, see e.g.~\cite[Section II.6]{tomDieck},
imply that $p$ may locally be inverted to obtain homeomorphisms to open subsets of the
individual sheets of the covering space. 
In particular, we may assume that the matrices $U_n$ converge in $\OO(d)$ to 
some $U\in H$ almost surely as~$n\to\infty$.

For each $a\in\Sd$ and $n\in\bN$ let $M_n^\pm(a):=\Med_\pm\bigl(P_n^a\bigr)$, 
and let $M(a)$ be the uniquely determined median of $P^a$. 
We define two stochastic processes $Y_n^\pm=\bigl(Y_n^\pm(a),a\in\Sd\bigr)$ with
index set $\Sd$ by
\begin{equation*}
  Y_n^\pm(a)=\sqrt{n}\big (M_n^\pm(a)-M(a)\big )\quad\text{for all } a\in \Sd,
\end{equation*}
and will use results presented by~\cite{Gr1996} 
and adopt arguments used by~\cite{KuelbsZinn} to prove that    
\begin{equation}\label{eq:FCLT}
        Y_n^\pm\dconv Y\quad \text{as }n\to\infty.
\end{equation}
Here $Y=(Y_a,a\in\Sd)$ is a centered Gaussian process with continuous
paths and covariance function
\begin{equation}\label{eq:FCLTcov}
     \cov(Y_a,Y_b) \, 
                  = \, \frac{K(a,b)}{f_a(M(a))f_b(M(b))},
\end{equation}
\vspace{-3mm}

\noindent
with
\vspace{-.5mm}
\begin{equation}\label{eq:defK}
   K(a,b)\, = \, P\bigl(\{x\in\bR^d:\, a\t x\le M(a), \, b\t x\le M(b)\}\bigr) 
                \,-\, \frac{1}{4}\,.
\end{equation}
We regard the $Y_n^\pm$ as processes with paths in
$\ell_\infty(\Sd)$, the space of real-valued bounded functions
on $\Sd$ endowed with the supremum norm, and the symbol '$\dconv$' refers to
weak or distributional convergence in the sense of Hoffmann-J{\o}rgensen; see
e.g.~\cite{vdV+W} and~\cite{Dudley} for details.   

For the proof of~\eqref{eq:FCLT} we start with a functional central limit theorem
for the empirical processes and then use an almost sure representation 
in order to be able to work with individual paths in an $\epsilon$-$\delta$ style.

Let $F(a,\cdot):=F_a(\cdot)$ be the distribution function of $a\t X$ and let
\begin{equation*}
   F_n(a,\cdot)\,=\,F_{a,n}(\cdot )\, =\, \frac{1}{n}\sum_{j=1}^n 1(a\t X_j\le \cdot\,)
\end{equation*}
be the (random) distribution function associated with $P_n^a$, $a\in\Sd$. Here $1(\,\cdot\,)$ 
denotes the indicator function of its (logical) argument.  We introduce the empirical process
$Z_n=\left (\sqrt{n}\left (F_n(a,t)-F(a,t)\right ),(a,t)\in \Sd\times \bR\right )$
as a process  with paths in $\ell_\infty(\Sd \times \bR)$. We obtain a semimetric $d$ on $\Sd\times \bR$ via
\begin{equation*}
 d\big ((a,s),(b,t)\big )^2
  = E\Bigl(\bigl(1(a\t X\!\le\! s) - P(a\t X\!\le\! s)\bigr)
         - \bigl(1(b\t X\!\le\! t) - P(b\t X\!\le\! t)\bigr)\Bigr)^2,
\end{equation*}
$(a,s),(b,t)\in \Sd\times \bR$. 
As noticed in \cite[Proof of Theorem~1]{Gr1996}, we then have $Z_n\dconv Z_0$, 
where $Z_0=(Z_0(a,t),(a,t)\in \Sd\times \bR)$ is a centered Gaussian process with 
covariance function
\begin{equation*}
  E\bigl(Z_0(a,s)Z_0(b,t)\bigr)\, =\, E\Bigl(\bigl(1(a\t X\!\le\! s) - P(a\t X\!\le\! s)\bigr)
         \bigl(1(b\t X\!\le\! t) - P(b\t X\!\le\! t)\bigr)\Bigr),
\end{equation*}
$(a,s),(b,t)\in \Sd\times \bR$, 
and sample paths that are bounded and
continuous with respect to $d$. The assumptions on $P$ ensure that the covariance 
function is continuous with respect to the usual topology on $\Sd\times \bR$. 
This implies that the process $Z_0$ has
continuous sample paths. 

For the remainder of the proof we abbreviate `almost surely' to `a.s.', 
and convergence refers to $n\to\infty$ unless specified otherwise. 

By the Skorokhod-Dudley-Wichura representation theorem~\cite[Theorem~3.5.1]{Dudley} 
there exists
a probability space $(\Tilde{\Omega},\Tilde{\cA},\Tilde{\bP})$ together with
perfect and measurable maps $g_n:\Tilde{\Omega}\rightarrow \Omega$, $n\in\bN_0$, 
such that $\Tilde{\bP}^{g_n}=\bP$ for all $n\in\bN_0$, and such that, with 
$\Tilde Z_n:= Z_n\circ g_n$,
\begin{equation}\label{eq:suptildeZ}
  \sup_{(a,t)\in\Sd\times\bR}\bigl|\Tilde{Z}_n(a,t)-\Tilde{Z}_0(a,t)\bigr| \, \rightarrow \, 0\quad \Tilde{\bP}\text{-{a.s.}}.
\end{equation}
As $F_n$ is a deterministic function of
$Z_n$ it follows that $\Tilde F_n:= F_n\circ g_n$ has the same distribution 
under $\Tilde \bP$ as $F_n$ under $\bP$, 
$n\in\bN$. This extends to the associated one-dimensional projections and their medians,
and to the $Y_n$-processes, where
$\Tilde{M}_n^\pm(a):=\Med_\pm\bigl(\Tilde F_n(a,\cdot)\bigr)$ and 
$\Tilde Y_n^\pm:=Y_n^\pm\circ g_n$, $n\in\bN$. 

Arguing as 
in~\cite[Proof of Theorem 1]{Gr1996} we obtain that 
\begin{equation*} 
     \tilde Y_n^\pm(a)\; \to\; 
     \tilde Y(a) := -f_a\bigl(M(a)\bigr)^{-1}\Tilde Z_0\bigl(a,M(a)\bigr)\quad 
                                                            \Tilde{\bP}\text{-a.s.}
\end{equation*}
for all $a\in\Sd$, and that 
\begin{equation}\label{eq:suptildeY}
  \sup_{a\in \Sd}\bigl|\Tilde{Y}_n^\pm(a)\bigr|\ = \ O(1)
  \quad \Tilde{\bP}\text{-{a.s.}}.
\end{equation}
For the integral of the marginal medians needed by~\cite{Gr1996} this was enough, but
here we need the  stronger stronger statement
\begin{equation}\label{eq:process_conv}
   \sup_{a\in \Sd}\bigl|\Tilde{Y}_n^\pm(a) - \Tilde{Y}(a)\bigr|\rightarrow 0\quad \Tilde{\bP}\text{-{a.s.}}.
\end{equation}
For the proof of~\eqref{eq:process_conv} we adopt some arguments given 
in~\cite[Proof of Proposition~1]{KuelbsZinn}. First we note that 
$F\bigl(a,M(a)\bigr)=1/2$ and then write
\begin{equation}\label{eq:entwicklung}
  F\bigl(a,M(a)+h\bigr) -  1/2 = f_a\bigl(M(a)\bigr)h + r(a,h)h,
\end{equation}
where $r(a,h)=f_a\bigl(M(a)+\zeta_{a,h}h\bigr)-f_a\bigl(M(a)\bigr)$ with 
some $\zeta_{a,h}\in [0,1]$. Using compactness of $\Sd$ we see that the
assumption~\eqref{as:cond2} implies 
\begin{equation*}
  \sup_{a\in\Sd}|r(a,h)|\le \sup_{a\in\Sd}\sup_{|s|\le |h|}\left |f_a\left (M(a)+s\right )-f_a\left (M(a)\right )\right |\rightarrow 0\hskip2mm \text{as}~h\to 0.
\end{equation*}
Further, \eqref{eq:suptildeY} yields
\begin{equation*}
  T_n\, := \, \sup_{a\in\Sd}\bigl|\Tilde{M}_n^\pm(a)-M(a)\bigr|\, \rightarrow\, 0\quad \Tilde{\bP}\text{-{a.s.}}
\end{equation*}
in view of the definition of the $Y$-processes.
Taken together this gives, $\Tilde{\bP}$-{a.s},
\begin{equation*}
 \sup_{a\in \Sd}\bigl|r\big (a,\tilde{M}_n^\pm(a)-M(a)\big )\bigr|\, 
  \le \sup_{a\in \Sd}\;\sup_{|s|\le T_n}\bigl|f_a\bigl(M(a)+s\bigr)-f_a\bigl(M(a)\bigr)
 \bigr|\, \rightarrow \, 0,
\end{equation*}
and using~\eqref{eq:suptildeY} again we get
\begin{align*}
  \Tilde{B}_n^\pm:=\sup_{a\in \Sd}|\Tilde{Y}_n^\pm(a)||r(a,\tilde{M}_n^\pm(a)-M(a))|\rightarrow 0\quad \Tilde{\bP}\text{-{a.s.}}.
\end{align*}
From~\eqref{eq:suptildeZ} and the continuity of 
$a\rightarrow \Tilde{Z}_0\left (a,M(a)\right )$ it follows that,
again $\Tilde \bP$-{a.s.},
\begin{equation}\label{eq:hilf1}
  \sup_{a\in \Sd}\left |\sqrt{n}\left (\Tilde{F}_{n}\big (a,\Tilde{M}_n^\pm(a)\big ) - F\big (a,\Tilde{M}_n^\pm(a)\big )\right ) -
      \Tilde{Z}_0\big (a,M(a)\big )\right |\; \rightarrow\;  0.
\end{equation}

Absolute continuity of $P$ implies that samples from $P$ are in general position
with probability~1, meaning that at most $d$ points are on a hyperplane. Hence
at most $d$ of the variables $a\t X_1,\ldots,a\t X_n$ coincide, and thus, given 
that the range of $\tilde F_n=F_n\circ g_n$ is contained in the range of $F_n$, 
\begin{equation*}
\bigl| \Tilde{F}_n\bigl(a,\Tilde{M}_n^\pm(a)\bigr)\, - \, 1/2\,\bigr| 
           \ \le \ d/n\,.
\end{equation*}
with probability~1. Using this and \eqref{eq:hilf1} we get
\begin{equation*}
  \sup_{a\in \Sd}\left |\sqrt{n}\Big (F\big (a,\Tilde{M}_n^\pm(a)\big ) - 1/2 \Big ) + \Tilde{Z}_0\big (a,M(a)\big )\right |\; \rightarrow\; 0 \quad \Tilde{\bP}\text{-a.s.}.
\end{equation*}
From
\begin{align*}
  \sqrt{n}\Bigl(F\big (a,\Tilde{M}_n^\pm(a)\big ) - 1/2 \Bigr) \,=\, \Tilde{Y}_n^\pm(a)\left (f_a\big (M(a)\big )+r\big (a,\Tilde{M}_n^\pm(a)-M(a)\big )\right ),
\end{align*}
which holds by \eqref{eq:entwicklung}, we obtain that, $\Tilde{\bP}$-{a.s.}, 
\begin{equation*}
\Tilde{A}_n^\pm:=\sup_{a\in \Sd} \left | \Tilde{Y}_n^\pm(a)\left (f_a\big (M(a)\big ) +
r\big (a,\Tilde{M}_n^\pm(a)-M(a)\big )\right ) + \Tilde{Z}_0\big (a,M(a)\big)\right |
   \rightarrow 0. 
\end{equation*}
Thus,
\begin{align*}
  \sup_{a\in\Sd}\left |\Tilde{Y}_n^\pm(a)f_a\big (M(a)\big )+\Tilde{Z}_0\big (a,M(a)\big )\right |
  \le \Tilde{A}_n^\pm+\Tilde{B}_n^\pm\rightarrow 0\quad \Tilde{\bP}\text{-{a.s.},} 
\end{align*}
which finishes the proof of~\eqref{eq:process_conv} as $a\mapsto f_a(M(a))$ is bounded
away from 0 on $\Sd$.

Returning to the un-tilded variables we see that this establishes the functional 
limit theorems in~\eqref{eq:FCLT}. 

We note in passing that  
$\sup_{a\in\Sd}\bigl|\Tilde Y_n^+(a)-\Tilde Y_n^-(a)\bigr|\,\rightarrow\, 0$ $\Tilde{\bP}$-{a.s.},
which implies
\begin{equation}\label{eq:shrink}
  \sup_{a\in\Sd}\left |\sqrt{n}\left (M_n^+(a) - M_n^-(a)\right ) \right |\rightarrow 0\quad
  \text{in }
    \bP\text{-probability}.
\end{equation}    
This shows that the volume of the interval $[\mMed_-(P_n^{U_n}),\mMed_+(P_n^{U_n})]$
shrinks fast enough so that the choice of a vector from this interval is irrelevant
for the distributional asymptotics.  

Now let $b_n(i)\t$ and $b(i)\t$ be the $i$th row of $U_n$, $n\in\bN$, and $U$
respectively, for $i=1,\dots,d$. In view of $U_n\to U$ in $\OO(d)$ 
we have $b_n(i)\to b(i)$ $\bP$-{a.s.} for $i=1,\ldots,d$. 
Let $W_n:=\sqrt{n}\eta_n = U_n\sqrt{n}(\theta_n-\theta)$ and 
\begin{align*}
  W_n^\pm:= \sqrt{n}  \begin{pmatrix} M_n^\pm(b_n(1))\\
                      \vdots\\
                      M_n^\pm(b_n(d))\end{pmatrix},
  \quad 
  W:=\begin{pmatrix} Y(b(1))\\ \vdots\\ Y(b(d)) \end{pmatrix}.  
\end{align*}
Note that $W_n$ is an element of the $d$-dimensional interval $[W_n^-,W_n^+]$ with 
probability~1. As the paths of the limit process are continuous the 
functional limit theorems~\eqref{eq:FCLT} now imply $W_n^\pm \dconv W$,
hence
\begin{equation}\label{eq:CLTspeziell}
         W_n \; \dconv\; W,
\end{equation}
with $W\sim N_d(0,\Xi)$, where the entries of $\Xi$ can be calculated 
from~\eqref{eq:FCLTcov} and~\eqref{eq:defK}. Indeed, as $M\equiv 0$ and as $(0,U)$
solves the quarter median problem for $P$, we obtain $K(b_i,b_i)=1/4$ and $K(b_i,b_j)=0$
if $i\not= j$, so that $\Xi=\Delta(P,U)$ with $\Delta(P,U)$ the diagonal matrix given
in~\eqref{eq:defDelta}.

We now apply $U_n\t$ on the left hand side, $U\t$ on the right
hand side of~\eqref{eq:CLTspeziell} and use~\cite[Theorem 5.5]{BillCoPM} again
to obtain the claimed asymptotic normality of $\sqrt{n}(\theta_n-\theta)$, together
with the mean and covariance matrix of the limit distribution.

It remains to show that $\Sigma(P,U)$ does not depend on the choice of $U$ from 
the set $H$. This follows easily from the fact that, for diagonal matrices $\Delta$,
we always have $P\t\Delta P=\Delta$ if $P$ corresponds to the permutation of two rows
or the multiplication by $-1$ of some row.

For the proof the assumption \eqref{as:cond4} is important. Without such an assumption we would need
the decomposition of $\sqrt{n}\bigl(M_n^\pm(b_n(i))) - M(b(i))\bigr)$ into two terms
$\sqrt{n}\bigl(M_n^\pm(b_n(i))) - M(b_n(i))\bigr)$ and 
$\sqrt{n}\bigl(M(b_n(i))) - M(b(i))\bigr)$. The first could still be analyzed 
with the above arguments, but for the second term the local behaviour of the function
$a\mapsto M(a)$  and distributional asymptotics of the estimates $U_n$ would become important;
see also Theorem~\ref{thm:alpha}. 

Further, \eqref{eq:shrink} and the ensuing comment show that the original estimate 
$\theta_n$ need not be permutation invariant. This may seem surprising 
as the proof relies on empirical process theory, which deals with the data
through their empirical distribution. Here it turned out to be enough that the estimate
can be squeezed between two other estimates that are such functions of the empirical 
distribution. 

\subsection{Proof of Theorem~\ref{thm:existell}}\label{subsec:proofexistell}
The transformation law for $d$-dimensional elliptical distributions, together with the 
equivariance from Proposition~\ref{prop:qumedequi}\,(a), means that we may assume that $\mu=0$ and that $\Sigma=\diag(\lambda_1,\ldots,\lambda_d)$. Recall that $e_1,\ldots,e_d$
is the canonical basis of $\bR^d$. For a subset $I$ of $\{1,\ldots,d\}$ let
\begin{equation*}
A_I := \bigl\{x\in\bR^d:\, x\t e_i\ge 0\text{ for }i\in I, x\t e_i\le 0\text{ for }i\notin I\bigr\}.
\end{equation*}
The assumptions imply that $P^T=P$ for all reflections $T:\bR^d\to\bR^d$
about hyperplanes $H(e_i)$,
hence $P(A_I)$ does not depend on $I$; also, $\sum_{I\subset\{1,\ldots,d\}}P(A_I)=1$.
A quarter space such as $\{x\in\bR^d:\, x\t e_i\ge 0, \, x\t e_j\ge 0\}$ with $i\not=j$  
can be written as the union of $2^{d-2}$ such sets $A_I$, where the intersections are all $P$-null 
sets, hence all of these have probability $2^{d-2}2^{-d}=1/4$. This shows that $(0,\Sigma)$
solves the quarter median problem for $P$.

In order to see that the quarter median is unique it is enough to show that $\CoMed(P^U)=0$
for all $U\in\OO(d)$ under the above assumptions on $P$. This follows
from the fact that $\Ell_d(0,\Sigma;h)$ is invariant under the reflection $x\mapsto -x$,
$x\in\bR^d$.

It remains to show that the solution of the quarter median is unique in the strictly elliptical 
case. To prove this, we recall that $\mu=0$ and that $\Sigma=\diag(\lambda_1,\ldots,\lambda_d)$.
Assume that $(0,[U])$ with $U\in \OO(d)$ is a solution 
of the quarter median problem. Let $b_1\t,\dots,b_d\t$ be the row vectors of $U$. 
Then $P^U=\Ell_d(0,\Tilde{\Sigma};h)$, where
$\Tilde{\Sigma}=\left (\Tilde{\sigma}_{ij}\right )_{1\le i,j\le d}=U\Sigma U^t$. 
With $U_{ij}$ as the $2\times d$ matrix with the rows $b_i\t,b_j\t,\,1\le i<j\le d,$ we see
that $P_{ij}:=P^{U_{ij}}$ is the elliptical distribution $\Ell_2(0,\Tilde{\Sigma}_{ij};h_2)$ 
with dispersion matrix 
\begin{equation*}
  \Tilde{\Sigma}_{ij}=U_{ij}\,\Sigma\,U_{ij}\t=
  \begin{pmatrix}\Tilde{\sigma}_{ii} & \Tilde{\sigma}_{ij}\\ \Tilde{\sigma}_{ji} & \Tilde{\sigma}_{jj}\end{pmatrix}
\end{equation*}
and density generator $h_2$ which is equal to $h$ if $d=2$ and given by
\begin{equation*}
        h_2(u)=\frac{\pi^{d/2-1}}{\Gamma(d/2-1)}\int_u^\infty
                      (t-u)^{d/2-2}\,h(t)\,dt\hskip2mm\text{for}~u>0,
\end{equation*}
if $d\ge 3$; see, e.g. \cite[Section 2.2.3]{Fang}. 
According to \cite[p.\,1466]{Gr1996} the probability assigned by $P_{ij}$ to the left lower quadrant is given by
\begin{equation*}
  \frac{1}{4}+\frac{1}{2\pi}\arcsin\left (\Tilde \sigma_{ij}/(\Tilde \sigma_{ii}\Tilde \sigma_{jj})^{1/2}\right ),
\end{equation*}
which is equal to $\frac{1}{4}$ if and only if $\Tilde\sigma_{ij}=0$, 
i.e.~if and only if $\Tilde{\Sigma}_{ij}$ is a diagonal matrix with positive diagonal elements. 
Therefore, $\Tilde{\Sigma}=U\Sigma U^t$ is a diagonal matrix,
$\Tilde{\Sigma}=\diag(\Tilde{\lambda}_1,\dots,\Tilde{\lambda}_d)$ say, or equivalently,
$\Sigma=U\t\Tilde{\Sigma}U=U\t\diag(\Tilde{\lambda}_1,\dots,\Tilde{\lambda}_d)U$.
It follows from this that the
set of eigenvalues $\bigl\{\Tilde{\lambda}_1,\ldots,\Tilde{\lambda}_d\bigr\}$ of
$\Tilde{\Sigma}$ and the set of eigenvalues $\bigl\{\lambda_1,\ldots,\lambda_d \bigr\}$ of
$\Sigma$ coincide and, as a consequence, that there is a permutation matrix $\Pi$ such that
$\Pi\t \Sigma \Pi=\Tilde{\Sigma}$, which in turn implies  
$U\t \Pi\t \Sigma U\Pi = \Sigma.$ Thus, putting $U\t\Pi\t=:Q=(q_{ij})$ we have
$Q\Sigma=\Sigma Q\t,$ i.e. $q_{ij}\lambda_j=\lambda_iq_{ij}.$ As the $\lambda_j$ are pairwise
distinct this gives $q_{ij}=0$ for $i\not=j$. Thus $Q$ is a diagonal matrix with the 
entries $\pm 1$ in the diagonal so that $U\sim I_d.$ Taken together
this shows that $(0,[I_d])$ is the unique solution of the quarter median problem. 

\subsection{Proof of Theorem~\ref{thm:CLTell}}\label{subsec:proofCLTell}
(a) We first assume that $\mu=0$ and that $\Sigma=\diag(\lambda_1,\ldots,\lambda_d)$.
Using Theorem~\ref{thm:CLTtheta} we only need to evaluate the diagonal elements 
of $\Delta(P,U)$, and because of Theorem~\ref{thm:existell} we may take $U$ to be the 
identity matrix. In particular, $b_i=e_i$ for $i=1,\ldots,d$.
Then $f_{e_i}$ is the $i$th marginal density of $P$, and the calculations preceding 
Theorem~\ref{thm:existell} lead to 
\begin{equation*}
     f_{e_i}\bigl(M(e_i)\bigr)\, =\, \frac{c_{d-1}(h)}{\sqrt{\lambda_i}c_d(h)},  
                       \quad i=1,\ldots,d.
\end{equation*}
Putting this together we see that, for such $P$, the asymptotic covariance matrix is 
indeed the specified multiple of $\Sigma$. 

We now use the equivariance properties of the quarter median to `bootstrap' this to 
general strictly elliptical distributions.
Hence suppose that $P=\Ell_d(\mu,\Sigma;h)$ and that 
$\lambda_1 > \lambda_2> \cdots> \lambda_d$ are the
eigenvalues of $\Sigma$. Then, for some $U\in\OO(d)$, 
$\Sigma=U\t\diag(\lambda_1,\ldots,\lambda_d)U$.
If $X_i$, $i\in\bN$, are independent with distribution $P$ then the random variables 
$Y_i:= U(X_i-\mu)$,  $i\in\bN$, are independent with distribution
$\Ell_d(0,\diag(\lambda_1,\ldots,\lambda_d);h)$. As 
$U(\theta_n - \mu\bigr)$ is a quarter median of $Y_1,\dots,Y_n$ the first part of 
the proof leads to 
\begin{equation*}
  \sqrt{n}\,U\bigl(\theta_n - \mu\bigr) \dconv N_d(0,\Xi),
\end{equation*}
with $\Xi=\sigma_{\QMed}^2\,\diag(\lambda_1,\ldots,\lambda_d)$.
Using the well-known properties of weak convergence and the behavior of multivariate 
normal distributions under affine transformations we finally obtain~\eqref{eq:CLT}.

\smallbreak
(b) Suppose that the eigenvalues of $\Sigma$
coincide and let $(U_n)_{n\in\bN}$ be as in the statement of the second part of the theorem. 
Then, for any subsequence $(U_{n'})$ there exists a subsubsequence $(U_{n''})$
such that $U_{n''}\to U_0$ for some 
$U_0\in\OO(d)$ as $n''\to\infty$. With $m$ for $n''$ and $b_m(i)\t$, $b(i)\t$
the $i$th row of $U_m$ respectively $U_0$, and using the same arguments as 
for~\eqref{eq:CLTspeziell}, we obtain
\begin{equation*}
   \sqrt{m} \begin{pmatrix} M_{m}(b_m(1))\\ \vdots\\
                            M_{m}(b_m(d))\end{pmatrix}
          \; \dconv\; \begin{pmatrix} Y(b(1))\\ \vdots\\ Y(b(d)) \end{pmatrix}
                        \quad\text{as }m\to\infty.
\end{equation*}
The symmetry of $P$ implies that the limit distribution does not depend on $U_0$.

\subsection{Proof of Theorem~\ref{thm:alpha}} \label{subsec:proofalpha}

As the result is somewhat tangential to the topic of multivariate location estimation
we only provide the main ideas together with some arguments specific to the quarter 
median application.

The general approach consists in writing the parameter $\theta=\theta(P)$ as the 
solution of an equation $\Psi_P(\theta)=z$ with a suitably chosen function 
$\Psi_P$ and known $z$. Let $\Psi_n$ be the empirical version of $\Psi_P$ and
suppose that $\Psi_n(\theta_n)$ sufficiently close to $z$ for large $n\in\bN$.
Then, under certain conditions, properties of the estimator sequence $(\theta_n)_{n\in\bN}$ can be 
obtained by localizing $\Psi_P$ near $\theta$ and then using the delta method. An excellent 
general exposition of this approach 
can be found in the textbook of~\cite{vdV}, see also~\cite{vdVSF}.
Specifically, \cite[Theorem 5.9]{vdV} together with \cite[Theorem 6.17]{vdVSF}
may serve as a guide towards filling in the details. Ultimately, this would also 
lead to an alternative proof for Theorem~\ref{thm:CLTell} in the case $d=2$
that does not make use of the results presented by~\cite{Gr1996}.

For $a,b>0$ let 
\begin{equation*}
     E(a,b)\, :=\, \Bigl\{(x,y)\t\in\bR^2:\, \frac{x^2}{a^2} +  \frac{y^2}{b^2} = 1\Bigr\}
             \, =\, \bigl\{(a\cos(t),b\sin(t))\t:\, 0\le t<2\pi\bigr\}
\end{equation*} 
be the boundary of the two-dimensional centered ellipse with main half-axes parallel to
the unit vectors $e_1$ and $e_2$ and of length $a$ and $b$ respectively.
On this set we consider the push-forward $Q(a,b):=\unif(\Seins)^{T_{a,b}}$ of the uniform distribution
on $\Seins=E(1,1)$ under the mapping $T_{a,b}:\bR^2\to\bR^2$, $(x,y)\t\mapsto (a x,b y)\t$. 
If $a\not= b$ then the solution of the quarter median problem for $Q(a,b)$ is unique and given by 
$\mu=(0,0)\t$ and $U=\diag(1,1)=U_\alpha$ with $\alpha=0$. 

Let $P$ be a probability measure on the Borel subsets of $\bR^2$. We introduce the function
$\Psi_P:\bR\times\bR\times [-\pi/4,\pi/4)\to \bR^3$, where the components
$\Psi_{P,i}$ are the respective probabilities of the upper right, upper left and lower left
quarter spaces associated with the argument. Formally, and with 
$b_\alpha,b_\alpha'$ the columns of $U_\alpha\t$, 
\begin{equation*}
   \Psi_{P,1}(x,y,\alpha) = P\bigl(V_{++}((x,y)\t;b_\alpha,b_\alpha')\bigr),
\end{equation*}   
and $\,\Psi_{P,2}(x,y,\alpha) = P\bigl(V_{-+}(\cdots)\bigr)$,  
$\,\Psi_{P,3}(x,y,\alpha) = P\bigl(V_{--}(\cdots)\bigr)$. 

Suppose that $P=Q(a,b)$ with $a\not=b$. Then elementary geometric considerations lead
to the following matrix of partial derivatives of $\,\Psi_{P}$ at $(x,y,z)=(0,0,0)$,
\begin{equation}\label{eq:Mmatrix}
M(a,b)\,  :=\,  \frac{1}{2\pi ab} 
  \begin{pmatrix} -b &b &b \\ -a & -a & a \\ b-a & a-b & b-a \end{pmatrix}.
\end{equation}       
For later use we  note that
\begin{equation*}
   M(a,b)^{-1} := \pi \begin{pmatrix} 0 & -b & ab/(b-a) \\ a & -b & 0 \\ a & 0 &ab/(b-a)\end{pmatrix}.
\end{equation*}  

In order to obtain the derivative for more general $P$  we use some structural properties 
of symmetric and elliptical distributions. Note that $Q(r,r)$ is the uniform distribution on 
the boundary of the  sphere with radius~$r$.  If $X\sim\Sym_2(h)$ then $\|X\|$ has density
\begin{equation*}
 f_{\|X\|}(r)\, =\, \frac{2\pi}{c_2(h)}\, r\,h(r^2),\quad r>0,
\end{equation*}
and given $\|X\|=r$, $X$ is uniformly distributed on $\{x\in\bR^2:\, \|x\|=r\}$.
Taken together this implies the mixture representation
\begin{equation*}
\Sym_2(h)\, =\, \frac{2\pi}{c_2(h)}\, \int_0^\infty Q(r,r) \, r \, h(r^2)\, dr.
\end{equation*}
Elliptical distributions are affine transformations of symmetric distributions.
In particular, 
\begin{equation*}
     P\, :=\,  \Ell_2\bigl((0,0)\t, \diag(a^2,b^2);h\bigr) \; =\; \Sym_2(h)^{T_{a,b}},
\end{equation*}
so that 
\begin{equation*}
        P \; =\; \frac{2\pi}{c_2(h)}\, \int_0^\infty Q(ar,br) \, r \, h(r^2)\, dr.
\end{equation*}   
Inserting quarter spaces we obtain a relation between $\Psi_P$ and the $\Psi$-functions 
corresponding to the measures $Q(a,b)$, $a,b>0$. This in turn can be used to relate the derivative 
$D(P)$ of $\Psi_P$ at $(0,0,0)$ to the 
derivatives obtained earlier for $Q(a,b)$, and we arrive at
\begin{equation*}
 D(P) \, =\, \frac{2\pi}{c_2(h)}\, \int_0^\infty M(ar,br)\, r \, h(r^2)\, dr
                              \, = \, \frac{\pi c_1(h)}{c_2(h)}\, M(a,b),
\end{equation*}
and thus
\begin{equation*}
 D(P)^{-1}  \, =\, \frac{c_2(h)}{c_1(h)}\, 
           \begin{pmatrix} 0 & -b & ab/(b-a) \\ a & -b & 0 \\ a & 0 &ab/(b-a)\end{pmatrix}.
\end{equation*}

If $(x,y,\alpha)$ solves the quarter median problem for an absolutely continuous $P$,
then $\Psi_P(x,y,\alpha)=(1/4,1/4,1/4)$.  
For a sample of size $n$ from $P$ the random vector of observation counts 
in the quarter spaces associated with the solution of the quarter median problem for $P$
has a multinomial distribution with parameters $n$
and $(1/4,1/4,1/4,1/4)$. The multivariate central limit theorem shows that the
standardized counts in the three quarter spaces corresponding to the components of $\Psi$
are asymptotically normal with mean vector 0 and covariance matrix  
\begin{equation*}
\Sigma_0 := \begin{pmatrix}
        3/16 & -1/16 &-1/16 \\ -1/16 & 3/16 & -1/16 \\ -1/16 & -1/16 & 3/16
\end{pmatrix}.
\end{equation*}  
We may now apply the delta method to obtain asymptotic normality for the triplet consisting 
of the quarter median coordinates and the rotation angle for a sequence of independent
random variables with distribution $P\,=\,\Ell_2\bigl((0,0)\t, \diag(a^2,b^2);h\bigr)$.
The limiting normal distribution is centered 
and has covariance matrix 
\begin{equation*}
\Sigma_1\, =\, (D(P)^{-1})\t \, \Sigma_0 \, D(P)^{-1} \; =\; \frac{c_2(h)^2}{4c_1(h)^2}\,
        \begin{pmatrix}
              a^2 & 0 & 0 \\ 0 & b^2 & 0 \\ 0 & 0 & a^2b^2/(a-b)^2
       \end{pmatrix}. 
\end{equation*}

It remains to extend this to $P=\Ell_2(\mu,\Sigma;h)$. We use essentially the same argument as 
at the end of Section~\ref{subsec:proofCLTell}\,(a). By assumption and Theorem~\ref{thm:existell}, we have that   
$\Sigma=U_\alpha\t\diag(\lambda_1,\lambda_2)U_\alpha$, where $U_\alpha$ is the
matrix representing the rotation by the angle $\alpha$ in the interior of
the half-open interval $I$, see \eqref{eq:parSO2}, and that $(\mu,U_\alpha)$
is the unique solution of the quarter median problem for $P$.  
Theorem~\ref{thm:consistency} implies that 
$\alpha_n-\alpha \in (-\pi/4,+\pi/4)$ with probability 1 for $n$ sufficiently large.
Let $Y_i:= U_\alpha\t(X_i-\mu)$, $i\in\bN$. Then, with probability 1 for $n$ large enough,  
$\left (U_\alpha\t(\QMed(X_1,\ldots,X_n)-\mu),U_{\alpha_n(X_1,\dots,X_n)-\alpha}\right )$ 
is a solution of the quarter median problem for $Y_1,\dots,Y_n$.  
Noting that $\lambda_1=a^2$ and $\lambda_2=b^2$ we may now apply the above asymptotic normality 
result to the $Y$-sequence to obtain
\begin{equation*}
  \begin{pmatrix}U_\alpha\t & 0 \\ 0 & 1\end{pmatrix}
  \begin{pmatrix}\sqrt{n}\big (\QMed(X_1,\ldots,X_n)-\mu\big )\\
    \sqrt{n}\big (\alpha_n(X_1,\ldots,X_n)-\alpha\big )\end{pmatrix}
  \ \dconv \ N_3(0,\Sigma_1).
\end{equation*}
Consequently,
\begin{equation}\label{final}
   \begin{pmatrix}\sqrt{n}\big (\QMed(X_1,\ldots,X_n)-\mu\big )\\
    \sqrt{n}\big (\alpha_n(X_1,\ldots,X_n)-\alpha\big )\end{pmatrix}
  \ \dconv \ N_3\left (0, \sigma_{\QMed}^2  \begin{pmatrix}\Sigma & 0 \\ 0 & \frac{\lambda_1\lambda_2}{(\sqrt{\lambda_1}
      -\sqrt{\lambda_2})^2}\end{pmatrix}\right ),
\end{equation}
which is the statement of the theorem.

\subsection{Proof of Proposition~\ref{prop:mess}}\label{subsec:proof_select}
We denote by $\mathcal C$ the set of non-empty closed subsets
of $\bR^d\times \SO(d).$ For $x_1,\ldots,x_n\in\bR^d$ let $\Gamma(x_1,\dots,x_n)$ 
be the (by Theorem~\ref{thm:main2} non-empty) set of solutions
$(\theta,U)\in \bR^d\times \SO(d)$ of the quarter
median problem for $P_{n;x_1,\dots,x_n}$. As permutations of $x_1,\ldots,x_n$
do not change $P_{n;x_1,\dots,x_n}$ we may regard $\Gamma$ as a function on
$E:=(\bR^d)^n/\sim$, where $(x_1,\ldots,x_n),(y_1,\ldots,y_n)\in (\bR^d)^n$ are
equivalent if
$y_i=x_{\pi(i)}$ for $i=1,\ldots,n$ 
for some permutation $\pi$ of the set $\{1,\ldots,n\}$. We aim to prove that the function
$\Gamma$ has values in $\mathcal C$ 
and that it is measurable in the sense that the set
\begin{equation}\label{eq:selcond}
    A(C)\, := \, \bigl\{ z\in E:\,  \Gamma(z)\cap C\neq\emptyset \bigr\}
\end{equation} 
is a Borel set for each $C\in \mathcal C$.
Then, by the measurable selection theorem of Kuratowski and Ryll-Nardzewski
(see, e.g. \cite[1.C Corollary]{Rockafellar}) there exists a measurable function
$\tau:E\to \bR^d\times \SO(d)$ such that $\tau(z)\in \Gamma(z)$ for
all $z\in E$. In fact, we will prove that the set
$A(C)$ is closed whenever $C\in \mathcal C$.

Arguing as in the first part of the proof of Theorem~\ref{thm:main2}
in Section~\ref{subsec:thm1} we see that if 
$(\theta_\ell,U_\ell)_{\ell\in\bN}$
is a sequence of elements in $\Gamma(z)$ converging to 
$(\theta,U)\in \bR^d\times \SO(d)$ as $\ell\to\infty$, then
$(\theta,U)\in \Gamma(z)$. This shows that the function $\Gamma$ has
values in $\mathcal C$. 

Now, more generally, fix $C\in \mathcal C$, 
let $A(C)$ be as in~\eqref{eq:selcond}, and let $(z_\ell)_{\ell\in\bN}$
be a sequence of elements of
$A(C)$ that converges to some element $z\in E$ as $\ell\to\infty$. 
As pointed out in connection with the transition from $\OO(d)$ to $\bH(d)$ 
in Section~\ref{subsec:proofCLTtheta} we can find 
$(x_1^{(\ell)},\dots,x_n^{(\ell)})\in z_\ell$, $l\in\bN$, and 
$(x_1,\dots,x_n)\in z$ such that $(x_1^{(\ell)},\dots,x_n^{(\ell)})$
converges to $(x_1,\dots,x_n)\in z$ in the space $(\bR^d)^n$. 
Then for each $\ell\in\bN$ there exists an associated element
$(\theta_\ell,U_\ell)$ of $\Gamma(x_1^{(\ell)},\dots,x_n^{(\ell)})$.
Arguing
as in Section~\ref{subsec:thm1} again, we see that along some subsequence 
$(\ell')$ it holds that $(\theta_{\ell'},U_{\ell'})\to (\theta,U)\in \bR^d\times \SO(d)$ 
and
\begin{equation*}
    P_{n;x_1^{(\ell')},\dots,x_n^{(\ell')}}^{U_{\ell'}}\ \rightarrow\  P_{n;x_1,\dots,x_n}^U
              \quad\text{weakly.}
\end{equation*} 
With $b_1\t,\dots,b_d\t$ as the row vectors of $U$ this gives  
$P_{n;x_1,\dots,x_n}\big (V_{\pm\pm}(\theta;b_i,b_j)\big )\ge \frac{1}{4}$. 
Consequently, $(\theta,U)\in \Gamma(z)$; also, $z\in A(C)$ as $C$ is closed.

\subsection{Proof of Theorem~\ref{thm:algo}}\label{subsec:proof_algo}
  Let $\theta$ be a quarter median of $P_n$, with associated orthogonal directions
  $b,b'\in S_2.$  Let $L_{\theta,b}=\{\theta + tb:t\in\bR\},\,L_{\theta,b'}=\{\theta + tb':t\in\bR\}$ 
  be the axes of the rectangular coordinate system with origin $\theta$ and directions $b,b'$.
  We consider two cases.

  (i) $\theta\in \{x_1,\dots,x_n\}$.
  We may assume without loss of generality that $\theta=x_1,$ and define
  \begin{equation*}
    \psi=\min\bigl(\min\bigl\{\arccos(|b_{1j}\t b|):2\le j\le n\bigr\},
          \,\min\bigl\{\arccos(|b_{1j}\t b'|):2\le j\le n\bigr\}\bigr).
  \end{equation*}
  There are orthogonal directions $b_*,b_*'\in \Seins$ obtained from $b,b'$ by the same clockwise or
  counterclockwise rotation with the angle $\psi$, where $b_*$ or $b_*'$ is an element of
  $\{-b_{12},b_{12},\ldots,-b_{1n},b_{1n}\}$. By construction, the number of data points $x_j$ lying in  
  $V_{\pm\pm}(\theta;b_{*},b_{*}')$ is greater than or equal to the number of data points $x_j$
  lying in $V_{\pm\pm}(\theta;b,b')$; also the number of data points $x_j$ lying in 
  $H_{\pm}(\theta;b_*)$ and the number of data points $x_j$ lying in 
  $H_{\pm}(\theta;b_*')$ is greater than or equal to the number of data points $x_j$
  lying in $H_{\pm}(\theta;b)$ and the number of data points $x_j$
  lying in $H_{\pm}(\theta;b')$, respectively. It follows that $(\theta,U_*)$ with
  $U_*$ as the matrix with the row vectors $b_*,b_*'$ solves the quarter median problem  
  for $P_n$. Thus by Proposition 3 (c) there is an $\eta\in \cM\bigl(P_n^{U_*}\bigr)$  
  such that $(\theta_*,U_*)$ with
  $\theta_*=U_*\t\eta$ also solves the quarter median problem for $P_n$. 
  Note that $\eta = \alpha b_* +\beta b_*'$, where $\alpha$ is median of $P_n^{b_*}$ 
  and $\beta$ is a  is median of $P_n^{b_*}$.
  Because $b_*$ or $b_*'$ is an
  element of $\{-b_{12},b_{12},\ldots,-b_{1n},b_{1n}\}$, and
  \begin{equation}\label{eq:reflect}
    \Med_\pm\left (P_n^{-b_*}\right )=-\Med_\mp\left (P_n^{b_*}\right ),
    \hskip5mm \Med_\pm\bigl(P_n^{-b_*'}\bigr)=-\Med_\mp\bigl(P_n^{b_*'}\bigr),
  \end{equation}
  it follows that $(\theta_*,U_*)\in \mathscr{L}(P_n).$

  (ii) $\theta\notin \{x_1,\dots,x_n\}$. With $\Tilde{\theta}$ as the
  intersection point of the line parallel to $L_{\theta,b}$ through the data point
  $\Tilde{x}\in \{x_1,\dots,x_n\}$ with the shortest distance to
  $L_{\theta,b}$, i.e. $\inf\{\|\tilde{x}-y\|:y\in L_{\theta,b}\}=\min\left \{\inf\{\|x_j-y\|:y\in L_{\theta,b}\}:1\le j\le n\right \}$,
  and the line parallel to $L_{\theta,b'}$ through the data point
  $\Tilde{\Tilde{x}}\in \{x_1,\dots,x_n\}$ with the shortest distance to
  $L_{\theta,b'}$, we have that $(\Tilde{\theta},U)$ with $U$ as matrix with the row vectors
  $b,b'$ is a solution of the quarter median problem for $P_n.$ If
  $\Tilde{x}=\Tilde{\Tilde{x}}$ then $\Tilde{\theta}=\Tilde{x}=\Tilde{\Tilde{x}}$
  and we arrive at case (i). So, we can (and do) assume that there are two points, $x_1,x_2$ say, such that
  $x_1\in L_{\theta,b},\,x_2\in L_{\theta,b'}.$
  Now we put $B:=\{b_{1j}:2\le j\le n\}\cup \{b_{2j}:3\le j\le n\}$ and define
  \begin{equation*}
    \psi=\min\bigl(\min\bigl\{\arccos(|d\t b|):d\in B\bigr\},
        \,\min\bigl\{\arccos(|d\t b'|):d\in B\bigr\}\bigr).
  \end{equation*}
  There are orthogonal directions $b_*,b_*'\in \Seins$ obtained from $b,b'$ by the same
  clockwise or
  counterclockwise rotation with the angle $\psi$, where $b_*$ or $b_*'$ is an element of $B\cup (-B)$.  
  With $L_1=\{x_1 + tb_*:t\in\bR\}$ and $L_2=\{x_2 + tb_*':t\in\bR\}$  
  we obtain the axes of a new rectangular coordinate system. The origin $\theta$ of
  the coordinate system with the axes $L_{\theta,b}$ and $L_{\theta,b'}$ and the origin
  $\theta_*$ of the new system are both located on the Thales circle about the line 
  segment connecting $x_1$ and $x_2$. 
  By construction, the number of data points $x_j$ lying in  
  $V_{\pm\pm}(\theta_*;b_*,b_*')$ is greater than or equal to the number of data points $x_j$
  lying in $V_{\pm\pm}(\theta;b,b')$. Further, the number of data points $x_j$ lying in 
  $H_{\pm}(\theta_*;b_*)$ and the number of data points $x_j$ lying in 
  $H_{\pm}(\theta_*;b_*')$ are greater than or equal to the number of data points $x_j$
  lying in $H_{\pm}(\theta;b)$ and the number of data points $x_j$
  lying in $H_{\pm}(\theta;b')$, respectively. Consequently, $(\theta_*,U_*)$ with
  $U_*$ as matrix with the row vectors $b_*,b_*'$ solves the quarter median for $P_n$.
  Again, by Proposition 3 (c), there is an $\eta\in \cM\bigl(P_n^{U_*}\bigr)$  
  such that $(\theta_{**},U_*)$, with $\theta_{**}=U_*\t\eta$ also
  solves the quarter median problem for $P_n$. Here, $b_*$ or $b_*'$ is an
  element of $B\cup (-B)$, so that, again by \eqref{eq:reflect}, we have
  $\bigl(\theta_{**},U_*\bigr) \in \mathscr{L}(P_n)$.

\vspace{.5cm}

\begin{acknowledgements} Dietrich Morgenstern (1924-2007) was the founding father of 
the Institut f{\"u}r Mathematische Stochastik at the then 
Technische Universit{\"a}t Hannover, which is now the Leibniz Universit{\"a}t Hannover,
and of two similar institutes in Freiburg and M{\"u}nster; 
see the obituary of~\cite{BGH2008}. 
Discussions with him were the starting point for the present paper.
\end{acknowledgements}

\section*{Conflict of interest}

The authors declare that they have no conflict of interest.


\bibliographystyle{spbasic}

\begin{thebibliography}{}

\bibitem[Bai et al(1990)]{Bai}
Bai ZD, Chen XR, Miao BQ, Rao CR (1990)
Asymptotic theory of least distances estimate in
multivariate linear models. Statistics 21:503--519

\bibitem[Baringhaus et al(2008)]{BGH2008}
Baringhaus L, Gr\"{u}bel R, Henze N (2008) Nachruf auf Dietrich Morgenstern. 
Jahresberichte der Deutschen Mathematiker-Vereinigung 110:101--113 

\bibitem[Billingsley(1968)]{BillCoPM}
  Billingsley P (1968) Convergence of probability measures. Wiley, New York

\bibitem[Brown(1983)]{BrownSMed}
  Brown BM (1983) Statistical uses of the spatial median. J Roy Statist Soc Ser B 45:25--30

\bibitem[Dhar and Chaudhuri(2011)]{DharChau}
  Dhar S, Chaudhuri P (2011) On the statistical efficiency of robust estimators of
  multivariate location. Stat Methodol 8:113--128 

\bibitem[Donoho and Gasko(1992)]{DG1992}
  Donoho DL, Gasko M (1992) Breakdown properties of location estimates based on
  halfspace depth and projected outlyingness. Ann Statist 20:1803--1827

\bibitem[Dudley(1999)]{Dudley}
  Dudley R (1999) Uniform central limit theorems. Cambridge University Press, Cambridge

\bibitem[Fang et al(1990)]{Fang}
  Fang TH, Kotz S, Ng KW (1990) Symmetric multivariate and related distributions.
  Chapman and Hall, London

\bibitem[G\'omez et al(1998)]{Gomez}
G\'omez E, G\'omez-Villegas MA, Mar\'in JM (1998) A multivariate generalization of the power
exponential family of distributions. Commun Statist - Theory Meth 27:589 -- 600

\bibitem[Gr\"ubel(1996)]{Gr1996}
Gr\"ubel R (1996) Orthogonalization of multivariate location estimators: the orthomedian. 
Ann Statist 24:1457--1473

\bibitem[Hall(2004)]{Hall}
Hall BC (2004) Lie groups, Lie algebras, and representations. Springer, New York

\bibitem[Kuelbs and Zinn(2013)]{KuelbsZinn}
Kuelbs J, Zinn J (2013) Empirical quantile CLTs for time-dependent data.
In: High dimensional probability VI:167--194, Progr Probab 66. Birkh{\"a}user/Springer,
Basel

\bibitem[Makeev(2007)]{Mak2007}
Makeev VV (2007) Equipartition of a continuous mass distribution.
J Math Sciences 140:551--557

\bibitem[Mass{\'e} and Plante(2003)]{MassPlant}
Mass\'e JC, Plante JF (2003) A Monte Carlo study of the accuracy and
robustness of ten bivariate location estimators.
Comput Stat Data Anal 42:1--26

\bibitem[Mass{\'e}(2004)]{Mass2004}
Mass\'e JC (2004) Asymptotics for the Tukey depth process, with an application to a multivariate 
trimmed mean.
Bernoulli 10:397--419

\bibitem[Mitchell(1989)]{Mitchell}
Mitchell AFS (1989) The information matrix, skewness tensor and a-connections for 
the general multivariate elliptic distribution.
Ann Inst Statist Math 41:289--304

\bibitem[Muirhead(1982)]{Muirhead}
Muirhead RJ (1982) Aspects of multivariate statistical theory. Wiley, New York  

\bibitem[Nolan(1999)]{Nolan}
  Nolan D (1999) On min-max majority and deepest points.
  Statist Probab Lett 43:325--333

\bibitem[Oja(1983)]{OjaMed}
  Oja H (1983) Descriptive statistics for multivariate distributions.
  Statist Probab Lett 1:327--333

\bibitem[Oja(2013)]{Oja}
  Oja H (2013) Multivariate median. In: Becker C, Fried R, Kuhnt S (Eds.)
  Robustness and complex data structures:\,3--15, Springer, Heidelberg

\bibitem[Oja and Niinima(1985)]{OjaNii}
  Oja H, Niinimaa A (1985) Asymptotic properties of the generalized
  median in the case of multivariate normality.
  J Roy Statist Soc Ser B 47:372--377

\bibitem[Rockafellar(1976)]{Rockafellar}
Rockafellar RT (1976) Integral functionals, normal integrands and measurable selections.
Nonlinear operators and the calculus of variations
(Summer School, Univ Libre Bruxelles, Brussels, 1975):157--207.
Lecture Notes in Mathematics 543, Springer, Berlin

\bibitem[Small(1990)]{Sm1990}
Small CG (1990) A survey of multidimensional medians. 
International Statistical Review 58:263--277

\bibitem[Somor\v{c}\'ik(2006)]{Somorcik}
Somor\v{c}\'ik J (2006) Tests using spatial median.
Austrian Journal of Statistics 35:331--338

\bibitem[tom Dieck(1991)]{tomDieck}
 tom Dieck T (1991) Topologie. Walter de Gruyter \& Co., Berlin

\bibitem[Tukey(1975)]{Tukey}
  Tukey JW (1975) Mathematics and the picturing of data. In: Proceedings of the 1974
  International Congress of Mathematics, Vancouver, Vol. 2, 523--531

\bibitem[van der Vaart(1998)]{vdV}
  van der Vaart AW (1998) Asymptotic statistics. Cambridge University Press, Cambridge

\bibitem[van der Vaart(2002)]{vdVSF}
  van der Vaart AW (2002) Semiparametric statistics. In: Bolthausen E, Perkins E,
  van der Vaart A. Lectures on probability theory and statistics. Lecture Notes in
  Mathematics 1781, Springer, Berlin

\bibitem[van der Vaart and Wellner(1996)]{vdV+W}
  van der Vaart AW, Wellner JA (1996) Weak convergence and empirical processes.
  Springer, New York
\end{thebibliography}

\end{document}